\documentclass[a4paper]{article}
\usepackage{amsmath,amssymb,amsfonts}
\usepackage{graphicx} 
\usepackage{bm}
\usepackage{subcaption}
\usepackage{placeins}
\usepackage{bbm}
\usepackage{cite}
\usepackage{enumerate}
\usepackage{hyperref}

\newcommand{\Z}{\bm{Z}}
\newcommand{\R}{\bm{R}}
\newcommand{\N}{\bm{N}}
\newcommand{\C}{\bm{C}}

\newcommand{\rmi}{{\rm i}}

\newenvironment{acknowledgements}{%
  \begin{abstract}
}{%
  \end{abstract}
}

\title{Accelerated Born series for high-frequency\\strong scattering by weakly non-radial media}

\author{Mirza Karamehmedovi\'c \thanks{Department of Applied Mathematics and Computer Science, Technical University of Denmark, Matematiktorvet 303B, DK-2800 Kgs. Lyngby, Denmark} \and Faouzi Triki \thanks{Laboratoire de Math\'ematiques et Applications, UMR CNRS 7348, Universit\'e
de Poitiers, TSA 61125, 11 Boulevard Marie et Pierre Curie, 86073 Poitiers
Cedex 9, France}}
\date{}

\begin{document}

\maketitle
\begin{abstract}
We present a practical method for accelerating the modified Born series for high-frequency Helmholtz scattering in one, two, and three dimensions. We consider bounded, nonnegative, compactly supported media, allowing small departures from radial symmetry in two and three dimensions. The method selects a single complex relaxation parameter using an approximation of the scattering operator's spectrum. Starting from precomputed eigenvalues for a constant reference medium, we estimate their evolution towards the prescribed radial profile. A leading-order Wentzel-Kramers-Brillouin (WKB) approximation reduces this spectral propagation to scalar nonlinear equations, followed by a first-order correction for the nonradial part of the medium. The estimated eigenvalues, together with a bound on the omitted spectrum, define a convex minimax problem for selecting the relaxation parameter. This spectral preparation is independent of the source and can therefore be reused for multiple scattering problems in the same medium. Numerical examples demonstrate convergence in strongly scattering cases where the scattering operator has spectral radius greater than one. Comparisons with direct integral-equation and finite-difference solvers assess accuracy, residual convergence, and computation time, and show computational savings at comparable accuracy in selected examples.

\end{abstract}

\section{Introduction}\label{sec:intro}

We here consider scattering problems in dimension $d=1,2,3$ at a fixed positive frequency (wavenumber) $k$. For the case $d=1$, fix a positive length $L$, a bounded and piecewise $C^1$ medium function $q$ supported in $[0,L]$, with
\[
q(x)\ge0,\quad\lim_{x\searrow0}q(x)\neq0,\quad\lim_{x\nearrow L}q(x)\neq0,
\]
as well as a source $f\in L^2(\R)$ supported in $[0,L]$. The scattering problem is then
\begin{equation}\label{eqn:PDE_d1}
\quad\left\{
\begin{array}{rcl}
u''+k^2(1+q(x))u&=&-f,\quad x\in (0,L),\\
-u'(0)&=&\rmi ku(0),\\
u'(L)&=&\rmi ku(L).
\end{array}
\right.
\end{equation}
For the case $d\in\{2,3\}$, let $B_R=\{x\in\R^d,\,\,|x|<R\}$ for some positive $R$, and fix a nonnegative-valued, bounded, and piecewise $C^1$ medium function $q$ supported in $\overline B_R$, with $\lim_{r\nearrow R}q(r\omega)\neq0$ for $\omega\in S^{d-1}$, as well as a source $f\in L^2(\R^d)$ supported in $\overline B_R$. Assume furthermore that $q$ is 'weakly non-radial', i.e., such that
\[
q(x)=q_{\rm rad}(|x|)+\eta Q(x),
\]
with the radial part
\[
q_{\rm rad}(r)=\frac{1}{|S^{d-1}|}\int_{\omega\in S^{d-1}}q(r\omega)d\omega,\quad r\in[0,R],
\]
and the perturbation $\|Q\|_{\infty}=1$, $|\eta|\ll1$. The scattering problem is then
\begin{equation}\label{eqn:PDE_d23}
\quad\left\{
\begin{array}{rcll}
\left[\Delta+k^2(1+q(x))\right]u&=&-f,\quad&x\in \R^d,\\
\lim_{|x|\rightarrow\infty}|x|^{(d-1)/2}(\partial_{|x|}u-\rmi ku)&=&0,\quad&\text{uniformly in}\,\,x/|x|\in S^{d-1}.
\end{array}
\right.
\end{equation}
In the following define the domain $D=(0,L)$ when $d=1$ and $D=B_R$ when $d\in\{2,3\}$. The problems~\eqref{eqn:PDE_d1} and~\eqref{eqn:PDE_d23} are posed in the usual outgoing weak sense, each PDE being supplied with the dimension-specific Sommerfeld radiation condition. The solution $u$ satisfies $u|_D\in L^2(D)$, and it is the total field resulting from the illumination of the medium $q$ by the source $f$.

\noindent
We are here concerned with the practical solvability of~\eqref{eqn:PDE_d1} and~\eqref{eqn:PDE_d23} in terms of the modified Born series (MBS) introduced by Kleinman \textit{et al.}~\cite{Kleinman-1988,Kleinman-1988-2,Kleinman-1990,Kleinman1990-2},
\begin{equation}\label{eqn:MBS}
u=\sum_{n=0}^{\infty}[M_{\gamma}(k)]^n\gamma V_1(k)f,
\end{equation}
where $M_{\gamma}(k)=(1-\gamma)I+\gamma k^2V_q(k)$, the constant $\gamma\in\C$ is the so-called preconditioner or overrelaxation parameter, and $V_1(k)$, $V_q(k)$ are the standard radiating volume potentials,
\[
V_a(k)u(x)=\int_D\Phi_k^{(d)}(x-y)a(y)u(y)dy,
\]
with
\[
\Phi_k^{(1)}(x)=(i/2k)e^{ik|x|},\quad\Phi_k^{(2)}(x)=\frac{i}{4}H_0^{(1)}(k|x|),\quad\text{and}\,\,\,\Phi_k^{(3)}(x)=\frac{e^{ik|x|}}{4\pi|x|},
\]
the unique outgoing fundamental solution of the Helmholtz operator $\Delta+k^2$ in $\R^d$ for $d=1,2,3$. Here $H_0^{(1)}$ is the Hankel function of the first kind and order zero.

\noindent
Convolving the PDE in \eqref{eqn:PDE_d1} or in~\eqref{eqn:PDE_d23} with the outgoing fundamental solution, and integrating by parts, we get the Lippmann-Schwinger integral equation
\begin{align} \label{intEq}
\left(I-k^2 V_q(k)\right) u = V_1(k)f, \qquad \textrm{in}\; D. 
\end{align}
It is well-known~\cite[Section 2]{2025-stochastic_spectral_decomp} that the integral equation \eqref{intEq}
is equivalent with \eqref{eqn:PDE_d1}-\eqref{eqn:PDE_d23}.  Recall~\cite{colton1998inverse} that the mapping $V_q(k): 
L^2(D) \to L^2(D)$  is compact, and that the Fredholm  equation of the second kind \eqref{intEq} has a unique solution $u = \left(I-k^2 V_q(k)\right)^{-1}V_1(k)f$.  If the spectral radius $\varrho(k^2V_q(k))$ of $k^2 V_q(k)$ is less than one then the solution of the integral equation is given by the Born series
\begin{equation*}
u \;=\; \sum_{n=0}^\infty   [k^2V_q(k)]^n(V_1(k)f).
\end{equation*}
In~\cite[Lemma 5]{Born_approximation} we proved, for scattering problem dimensions $d=1,2,3$ and general bounded media $q$, the existence of $\gamma$ that ensures the convergence of the modified Born series~\eqref{eqn:MBS} regardless of the value of the spectral radius $\varrho(k^2V_q(k))$, that is, including the strong scattering case (large $\|q\|_{\infty}$) and the high-frequency case (large $k$). Moreover, in~\cite[Lemma 6]{Born_approximation} we provided a practically calculable sufficient condition for such $\gamma$ in dimension $d=1$ with $q(x)=\|q\|_{\infty}\bm 1_{[0,L]}(x)$.\\

Ideally, we estimate the optimal preconditioner parameter $\gamma_{\rm opt}$ by minimizing the spectral radius of $M_{\gamma}(k)$; in practice, we solve the minimax problem~\eqref{eq:gamma_opt_definition} using the estimated spectrum. The aim is to improve the worst-case asymptotic convergence rate of the MBS~\eqref{eqn:MBS} and thereby reduce the number of terms needed to reach a prescribed tolerance. To this end, we estimate in Section~2 the spectral evolution of the
operator $k^2V_{q_t}(k)$ for the homotopy
\[
q_t=(1-t)\|q\|_\infty\mathbf 1_D+tq_{\rm rad},\qquad t\in[0,1],
\]
followed by a correction, to first order in $\eta$, of the computed eigenvalues $\lambda_j(1)$ from those of $k^2V_{q_{\rm rad}}(k)$ to those of $k^2V_q(k)$; recall that we assume $q(x)=q_{\rm rad}(|x|)+\eta Q(x)$ with $\|Q\|_{\infty}=1$ and with small $\eta$. To leading order in a WKB approximation~\cite{Babich1999}, the spectral evolution of $k^2V_{q_t}(k)$ from $t=0$ to $t=1$ can be
computed by solving a collection of scalar nonlinear
equations for the endpoint eigenvalues.
In Section~\ref{sec:estimating}, we use the estimated spectrum of $k^2V_q(k)$ to arrive at inexpensive numerical estimates of $\gamma_{\rm opt}$. In Section~\ref{sec:numerical}, we present numerical examples of the application of our accelerated MBS solver, including validation against direct solvers.\\

\noindent
Our proposed workflow is thus as follows:
\begin{enumerate}
    \item Pre-compute eigenvalues $\lambda_j(0)$, $j=1,\dots,J$, of $k^2V_{q_0}(k)$ for the values $k$, $\|q\|_{\infty}$, and $|D|$ of interest.
    \item Given a specific medium $q=q_{\rm rad}+\eta Q$ of interest, propagate the selected
eigenvalue branches for $k^2V_{q_t}(k)$, with $q_t=(1-t)q_0+tq_{\rm rad}$, from $t=0$ to $t=1$ using the high-frequency
spectral-evolution equation~\eqref{eqn:ode}. This reduces
to the solution of $J$ scalar nonlinear equations, as specified in~\eqref{eqn:use0} for $d=1$, in~\eqref{eqn:use1} for $d=2$, and in~\eqref{eqn:use2} for $d=3$.
\item If the medium $q$ is nonradial ($\eta\neq0$) then correct the computed $\lambda_j(1)$, $j=1,\dots,J$, using~\eqref{eqn:mu}.
    \item Use the (corrected) $\lambda_j(1)$, $j=1,\dots,J$, to estimate $\gamma_{\rm opt}$.
    \item For each error tolerance $\varepsilon$ and each specific source $f$ of interest, compute
    \[
    u_{\rm MBS}=\sum_{n=0}^{N(\varepsilon)}[M_{\gamma_{\rm opt}}(k)]^n\gamma_{\rm opt}V_1(k)f.
    \]
\end{enumerate}
Step 1 encompasses a large class of media, as well as all sources. Steps 2--4 encompass all sources for a given medium.  Thus, if we for example consider scattering by a fixed medium excited by $\nu$ different sources then the required computations reduce to a small fixed overhead (steps 1--4), followed by $\nu\times(N(\varepsilon)+1)$ matrix-vector multiplications and $\nu\times N(\varepsilon)$ vector summations. The most computationally intensive part is thus well-suited for GPU implementation.

\noindent
The systems~\eqref{eqn:PDE_d1}-\eqref{eqn:PDE_d23} arise in the time-harmonic case, with the time-dependence factor $e^{-\rmi \omega t}$, of the initial-value problem for the scalar wave equation
\begin{equation*}
\left\{\begin{array}{rcl}u_{tt}-c^2(x)\Delta u&=&0,\\
u(x,0)&=&0, \\ \partial_t u(x, 0)&=&c^2 f,\end{array}\right.
\end{equation*}
where the medium function is given by $q=c_0^2/c^2-1$, 
with $c_0>0$ being the ambient constant wavespeed, $c(x)>0$ the variable wavespeed describing the scatterer, and $k=\omega/c_0$.

\par
\noindent Kleinman \textit{et al.}~\cite{Kleinman-1988,Kleinman-1988-2,Kleinman-1990,Kleinman1990-2} devised the modified Born series by introducing an overrelaxation of the original Lippmann-Schwinger equation $(I-k^2V_q(k))u=\psi$ in terms of the equivalent equation $(I-M_{\gamma}(k))u=\gamma\psi$ with $M_{\gamma}(k)=(1-\gamma)I+\gamma k^2V_q(k)$, and where $\gamma\in\C$ is the overrelaxation parameter (see also~\cite{Engl-Nashed-1981}). For scattering configurations in dimensions $d=1,2,3$ and involving piecewise H\"older continuous contrasts, they proved the existence of complex constants $\gamma$ for which the modified Born series $u=\sum_{n=0}^{\infty}[M_{\gamma}(k)]^n\gamma\psi$ converges, as well as a sufficient condition for the choice of overrelaxation parameter resulting in convergent modified Born series. However, this sufficient condition involves parameters defining the geometry of the spectrum of the operator $k^2V_q(k)$, making the result less readily applicable in practice. Kleinman \textit{et al.}~\cite{Kleinman-1990} therefore also supplied an empirical recipe for the choice of the overrelaxation parameter, as well as some supporting numerical evidence. More recently, Osnabrugge et al.~\cite{Osnabrugge-2016} devised a more elaborate approach, with a non-constant function ('preconditioner' $\gamma(x)$) taking the place of the overrelaxation parameter, and using operators of the form $F^{-1}\widetilde{g}_0F$, with $F,F^{-1}$ the Fourier transform and its inverse, respectively, and $\widetilde g_0$ the Fourier transform of the outgoing fundamental solution of the Helmholtz operator in $\R^d$ at a complex frequency $k_0^2+i\varepsilon$. The added numerical complexity of this approach pays off in that the function $\gamma(x)$ is well-known from the outset and it is therefore not required to estimate a good overrelaxation parameter. Finally, Stanziola \textit{et al.}~\cite{Arridge-2023} used machine learning to arrive at Learned Born Series (LBS), speeding up the convergence of the MBS from~\cite{Osnabrugge-2016} for strongly scattering media.

\section{WKB approximation of the spectral evolution}\label{sec:WKB}
In this section we describe the two distinct steps outlined in Section~\ref{sec:intro}, namely the spectral propagation of $k^2V_{q_t}(k)$ for $q_t=(1-t)q_0+tq_{\rm rad}$, $t\in[0,1]$, and the subsequent first-order correction of the spectrum of $k^2V_{q_{\rm rad}}(k)$ to that of $k^2V_q(k)$, where $q(x)=q_{\rm rad}(|x|)+\eta Q(x)$ with $q_{\rm rad}$ radial, $\|Q\|_{\infty}=1$, and $\eta$ small.

\subsection{Spectral propagation}\label{sec:sp}

Given a medium $q$ as described in Section~\ref{sec:intro}, write $q_0=\|q\|_{\infty}\bm1_{\overline D}$ and consider the homotopy
\begin{equation}\label{eqn:homotopy}
q_t = (1-t)q_0+tq_{\rm rad}, \qquad t\in[0,1].
\end{equation}
Note that, for all $t\in[0,1]$, we have $q_t(x)\ge0$ and $q_t$ is nowhere vanishing at $\partial D$. We are interested in the one-parameter family of eigenvalue problems
\begin{equation}\label{eqn:1pevp}
k^2V_{q_t}(k)\varphi_t=\lambda(t)\varphi_t,\quad t\in[0,1].
\end{equation}
Specifically, we seek the variation of eigenvalue branches $\lambda(t)$ along the homotopy $q_t$, $t\in[0,1]$. Since $V_{q_t}(k) = V_{q_0}(k)+ tV_{q_{\rm rad}-q_0}(k)$ is $t$-analytic, $\lambda(t)$ is piecewise analytic, with at most finitely many exceptional points \cite{kato1966perturbation,ammari2004splitting}. In dimension $d=1$ the eigenvalues are geometrically simple because of the Robin boundary condition at $x=0$. Indeed, assume $\varphi_t^{(j)}$, $j=1,2$, are two eigenfunctions corresponding to the same eigenvalue $\lambda(t)$, with $\varphi_t^{(j)}(0)=a_j$. Then $a_1\neq0$ and $a_2\neq0$, and
\[
\varphi_t^{(1)}(0)=\frac{a_1}{a_2}\varphi_t^{(2)}(0)\quad\text{and}\quad\left(\varphi_t^{(1)}\right)^\prime(0)=\frac{a_1}{a_2}\left(\varphi_t^{(2)}\right)^\prime(0),
\]
so that $\varphi_t^{(1)}=(a_1/a_2)\varphi_t^{(2)}$. In dimensions $d=2,3$, the reference medium $q_0$ is radial, and the corresponding Helmholtz operator $\Delta+k^2(1+q_0/\lambda(0))$ is rotationally invariant. Hence its eigenspaces decompose into angular Fourier modes in \(d=2\) and spherical harmonics in \(d=3\). For $d=2$, the modes are labeled by $m\in\Z$, and for $d=3$ they are labeled by the angular momentum $\ell\in\N_0$. For $d=2$ the eigenvalues have multiplicity 1 when $m=0$ and 2 when $m\neq0$, and for $d=3$ each corresponding eigenvalue carries the multiplicity $2\ell+1$. Since \(q_{\rm rad}\) is radial, these multiplicities are preserved along the homotopy.
\noindent
Taking variations of~\eqref{eqn:1pevp}, we get
\[
k^2V_{q_t+\delta q}(k)(\varphi_t+\delta\varphi)=(\lambda(t)+\delta\lambda)(\varphi_t+\delta\varphi),
\]
so, to leading order, we have
\[
k^2V_{q_t}(k)\delta\varphi+k^2V_{\delta q}(k)\varphi_t=\lambda(t)\delta\varphi+\delta\lambda\varphi_t.
\]
Now use the bilinear pairing $\langle u,v\rangle=\int_Du(x) {v(x)}dx$ and consider the transpose $(k^2V_{q_t}(k))^T$, defined by $\langle k^2V_{q_t}(k)u,v\rangle=\langle u,(k^2V_{q_t}(k))^Tv\rangle$. We readily check that the eigenvalues and eigenvectors of $(k^2V_{q_t}(k))^T$ are given by $\lambda(t)$ and $q_t\varphi_t$, respectively: writing the multiplication operator $Q_{q_t}=q_t\cdot(\cdot)$, and $K(k)=k^2\int_D\Phi_k^{(d)}(x-y)(\cdot)(y)dy$, we have $Q_{q_t}^T=Q_{q_t}$ and $K(k)^T=K(k)$, so $k^2V_{q_t}(k)\varphi_t=\lambda(t)\varphi_t$ implies
\[
(k^2V_{q_t}(k))^T(q_t\varphi_t)=Q_{q_t}K(k)(q_t\varphi_t)=Q_{q_t}k^2V_{q_t}(k)\varphi_t=\lambda(t)q_t\varphi_t.
\]
Furthermore, we have $\delta q=(q_{\rm rad}-q_0)\delta t$ for every $t\in(0,1)$, so
\[
\delta t\langle k^2V_{q_{\rm rad}-q_0}(k)\varphi_t,q_t\varphi_t\rangle=\langle k^2V_{\delta q}(k)\varphi_t,q_t\varphi_t\rangle=\delta\lambda\langle\varphi_t,q_t\varphi_t\rangle,
\]
which leads to
\[
\lambda'(t)=\frac{\langle k^2V_{q_{\rm rad}-q_0}(k)\varphi_t,q_t\varphi_t\rangle}{\langle\varphi_t,q_t\varphi_t\rangle}=\lambda(t)\frac{\int_D(q_{\rm rad}(x)-q_0(x))\varphi_t^2(x)dx}{\int_Dq_t(x)\varphi^2_t(x)dx},\quad t\in(0,1),
\]
or equivalently
\begin{equation}\label{eqn:ode}
\frac{d}{dt}\log\lambda(t)=\frac{\int_D(q_{\rm rad}(x)-q_0(x))\varphi_t^2(x)dx}{\int_Dq_t(x)\varphi^2_t(x)dx},\quad t\in(0,1).
\end{equation}

\noindent
We do not have a formal justification for the division by $\int_Dq_t\varphi_t^2dx$; instead, we shall assume this division is well-defined, and rely on numerical examples to check the practical applicability of the approach. The idea now is to start with the eigenvalues of $k^2V_{q_0}(k)$, at $t=0$, which is a canonical case, and then to propagate these eigenvalues to those of $k^2V_{q_{\rm rad}}(k)$, at $t=1$, using the ODE~\eqref{eqn:ode}. This of course requires some knowledge of the point spectrum of $k^2V_{q_0}(k)$, as well as of the eigenfunctions $\varphi_t$ for $t\in[0,1]$. In what follows, we shall approximate the eigenfunctions $\varphi_t$ using the high-frequency Wentzel–Kramers–Brillouin (WKB,~\cite{Babich1999}) formalism. Regarding the eigenvalues of $k^2V_{q_0}(k)$, a finite representative number of these can be computed numerically once and for all for the contrast $\|q\|_{\infty}\bm1_{\overline D}$ and the size parameter $\kappa=kL$ ($d=1$) or $\kappa=kR$ ($d=2,3$) of interest, or they can be estimated semi-analytically. In dimension $d=1$, for example, we already have an explicit parametrization~\cite[proof of Proposition 3]{Born_approximation} of the curve $\Lambda_0$ on which the eigenvalues $\lambda_j(0)$ of $k^2V_{q_0}(k)$ reside, namely
\[
\Lambda_0=\left\{\frac{\|q\|_{\infty}}{s^2-1},\,\,s\in S_-\cup S_+\right\},
\]
with
\[
S_{\pm}=\left\{\frac{-\cosh\kappa t\pm\sqrt{1-t^2\sinh^2\kappa t}}{\sinh\kappa t}+it,\,\,t\in(0,T]\right\},
\]
and $T\sinh(\kappa T)=1$. In dimensions $d=2$ and $d=3$, we could attempt to find curves in $\C$ that envelop the point spectrum $\sigma_p(k^2V_{q_0}(k))$ by solving transcendental equations involving (spherical) Bessel and Hankel functions. Instead of this approach, however, we here opt for the one-off numerical computation of a finite number of eigenvalues of $k^2V_{q_0}(k)$ for all $d\in\{1,2,3\}$.

\noindent
We do not in general have an analytic expression of $\varphi_t$ for $t\in(0,1]$. To continue the analysis of~\eqref{eqn:ode}, we first note that the eigenvalue problem~\eqref{eqn:1pevp} is equivalent with the variable-coefficient Helmholtz system
\begin{equation}\label{eqn:system}
\left\{\begin{array}{c}\left[\Delta+k^2(1+q_t(x)/\lambda(t))\right]\varphi_t=0,\qquad x\in D,\\\text{+ Sommerfeld radiation condition},\end{array}\right.
\end{equation}
which has a nontrivial solution only if $\Im\lambda(t)>0$, that is, only in the case of an exponentially unstable (active medium) wave equation -- see, e.g., the proof of Lemma 5 in~\cite{Born_approximation}. Then, in what follows, we will use the high-frequency ($k\gg1$) WKB ansatz~\cite{Babich1999} to arrive at approximate solutions of~\eqref{eqn:system}.
\subsubsection{Case $d=1$}\label{subsec:1D}
Now first let $d=1$ and \(q=q_{\rm rad}\in C^1([0,L])\). For piecewise \(C^1\) contrasts, the WKB construction is understood on each smooth subinterval. Interface matching terms may generate additional reflected and transmitted contributions; these are here neglected at the leading averaged level. In one dimension we have only two {\it rays}, and the system~\eqref{eqn:system} is readily solvable to leading order in $k$, giving
\[
\varphi_t^2(x)\approx\frac{c_1^2(t)e^{2ikS_{t,1}(x)}+c_2^2(t)e^{-2ikS_{t,1}(x)}+2c_1(t)c_2(t)}{(1+q_t(x)/\lambda(t))^{1/2}},
\]
with
\[
S_{t,j}(x)=(-1)^{j+1}\int_0^x(1+q_t(y)/\lambda(t))^{1/2}dy,\,\,\,a_{t,j}(x)=\frac{c_j(t)}{(1+q_t(x)/\lambda(t))^{1/4}},\,\,\,j=1,2.
\]
Assuming that the WKB phases have no stationary points on \((0,L)\), that
the relevant branch of \(S_{t,1}\) has controlled imaginary part, and that
the oscillatory terms are negligible after high-frequency averaging, we obtain
formally, after substitution into~\eqref{eqn:ode}, and to leading order as $k\rightarrow\infty$, that
\[
\frac{\lambda'(t)}{\lambda(t)}
=
\frac{
\displaystyle\int_0^L
\frac{q_{\rm rad}(x)-q_0(x)}
{\left(1+q_t(x)/\lambda(t)\right)^{1/2}}\,dx
}{
\displaystyle\int_0^L
\frac{q_t(x)}
{\left(1+q_t(x)/\lambda(t)\right)^{1/2}}\,dx
}.
\]
We here also assume that both WKB components are present, so that
\(c_1(t)c_2(t)\neq 0\) along the branch under consideration. Throughout, the square root branch is chosen continuously along the eigenvalue branch $\lambda(t)$, with the initial branch fixed at $\lambda(0)$.
Now with
\[
H(t)=\int_0^L\left(1+\frac{q_t(x)}{\lambda(t)}\right)^{1/2}dx,
\]
we have, to leading order as $k\rightarrow\infty$,
\begin{align*}
\frac{dH}{dt}(t)&=\frac{1}{2\lambda(t)}\int_0^L\frac{q_{\rm rad}(x)-q_0(x)}{(1+q_t(x)/\lambda(t))^{1/2}}dx-\frac{\lambda'(t)}{2\lambda(t)^2}\int_0^L\frac{q_t(x)}{(1+q_t(x)/\lambda(t))^{1/2}}dx\\&=0,
\end{align*}
so specifically
\begin{equation}\label{eqn:use0}
\int_0^L\left(1+\frac{q_{\rm rad}(x)}{\lambda(1)}\right)^{1/2}dx-L\left(1+\frac{\|q\|_{\infty}}{\lambda(0)}\right)^{1/2}=H(1)-H(0)=o(1),\,\,\,k\rightarrow\infty,
\end{equation}
which can be solved for $\lambda(1)$ using a suitable zero-finding algorithm.

\subsubsection{Case $d=2$}\label{subsec:d=2}

Using the polar coordinates \((r,\theta)\), separate the variables in~\eqref{eqn:system} to get
\begin{equation}\label{eqn:phi}
\varphi_t(r,\theta)=\sum_{m=0}^{\infty}r^{-1/2}v_{m,t}(r)(a_{m,t}\cos m\theta+b_{m,t}\sin m\theta),
\end{equation}
where
\[
v_{m,t}''(r)+k^2\left[1+\frac{q_t(r)}{\lambda(t)}-\frac{m^2-1/4}{k^2r^2}\right]v_{m,t}(r)=0,\qquad 0<r<R,\,\,m\in\N_0,
\]
with the regularity condition at the origin given by
\[
r^{-1/2}v_{m,t}(r)=O(r^m),\qquad r\searrow0,
\]
and the outgoing condition for $r\rightarrow\infty$, as well as the continuity of $v_{m,t}$ and $v'_{m,t}$ at \(r=R\), implying
\[
v_{m,t}'(R)=\left(\frac{1}{2R}+k\frac{(H_m^{(1)})'(kR)}{H_m^{(1)}(kR)}\right)v_{m,t}(R),\qquad m\in\N_0.
\]
Here \(H_m^{(1)}\) is the Hankel function of the first kind and order
\(m\). As long as we do not encounter exceptional points, which we in the following assume is the case, we can describe the evolution of $\lambda(t)$ using the spectral propagation formula~\eqref{eqn:ode} with any of the eigenvectors $\varphi_{m,t}(r,\theta)=r^{-1/2}v_{m,t}(r)(a_m\cos m\theta+b_m\sin m\theta)$ contributing nontrivially to~\eqref{eqn:phi}. Since we wish to use~\eqref{eqn:ode}, we also assume that, for the chosen index $m$,
\[
\int_Dq_t(r)\varphi_{m,t}(r,\theta)^2rdrd\theta\neq0.
\]
Thus, with
\[
p_{m,t}(r;z)=\left(1+\frac{q_t(r)}{z}-\frac{m^2-1/4}{k^2r^2}\right)^{1/2}
\]
and
\[
S_{m,t}(r)=\int_{r_{\ast}}^rp_{m,t}(\varrho;\lambda(t))\,d\varrho,
\]
for an appropriate positive $r_{\ast}$, the WKB approximation is to leading order given by
\[
v_{m,t}(r)\approx p_{m,t}(r;\lambda(t))^{-1/2}\left(c_m(t)e^{ikS_{m,t}(r)}+        d_m(t)e^{-ikS_{m,t}(r)}\right).
\]
Substitution into the spectral evolution equation \((2.3)\), followed by
the removal of the oscillatory terms as in the one-dimensional case, and the  cancellation of the cross terms due to angular orthogonality, gives
\[
\frac{\lambda'(t)}{\lambda(t)}=\frac{\displaystyle\int_0^R\frac{q_{\rm rad}(r)-q_0}{p_{m,t}(r;\lambda(t))}\,dr}{\displaystyle\int_0^R\frac{q_t(r)}{p_{m,t}(r;\lambda(t))}\,dr}.
\]
Define
\[
H_m(t)=\int_0^R\left(p_{m,t}(r;\lambda(t))-p_{m,0}(r;\lambda(0))\right)dr.
\]
Then, to leading order as $k\rightarrow\infty$,
\[
\frac{dH_m}{dt}(t)=\frac{1}{2\lambda(t)}\int_0^R\frac{q_{\rm rad}(r)-q_0}{p_{m,t}(r;\lambda(t))}dr-\frac{\lambda'(t)}{2\lambda(t)^2}\int_0^R\frac{q_t(r)}{p_{m,t}(r;\lambda(t))}dr=0,
\]
and therefore
\begin{align}\label{eqn:use1}
\int_0^R(p_{m,1}\left(r;\lambda(1))-p_{m,0}(r;\lambda(0))\right)dr=H_m(1)=o(1),\quad k\rightarrow\infty.
\end{align}
This equation can be solved numerically for \(\lambda(1)\), using a zero finding algorithm and the
precomputed eigenvalue \(\lambda(0)\) of \(k^2V_{q_0}(k)\).

\subsubsection{Case $d=3$} Here we have a development and argumentation analogous to the case $d=2$: with $(r,\theta,\phi)$ the spherical coordinates, and assuming no exceptional points are encountered, we can track $\lambda(t)$ using~\eqref{eqn:ode} with any appropriate eigenvector
\[
\varphi_{\ell,m,t}(r,\theta,\phi)=r^{-1}v_{\ell,t}(r)Y_{\ell}^m(\theta,\phi),\quad \ell\in\N_0,\,\,m\in\{-\ell,\dots,\ell\},
\]
where
\[
v_{\ell,t}''(r)+\left(k^2(1+q_t(r)/\lambda(t))-\ell(\ell+1)/r^2\right)v_{\ell,t}(r)=0,\quad 0<r<R,
\]
and $v_{\ell,t}(r)=O(r^{\ell+1})$ for $r\searrow0$, as well as where
\[
v_{\ell,t}'(R)=\left(\frac{1}{R}+k\frac{(h_{\ell}^{(1)})'(kR)}{h_{\ell}^{(1)}(kR)}\right)v_{\ell,t}(R),
\]
and
\[
\int_Dq_t(r)\varphi_{m,\ell,t}(r,\theta,\phi)^2r^2\sin\theta drd\theta d\phi\neq0.
\]
Here $Y_{\ell}^m$ is a real spherical harmonic, and $h_{\ell}^{(1)}$ is the spherical Hankel function of order $\ell$ and of the first kind. The WKB approximation is then, to leading order,
\[
v_{\ell,t}(r)=p_{\ell,t}(r;\lambda(t))^{-1/2}\left(c_{\ell}(t)e^{ikS_{t,\ell}(r)}+d_{\ell}(t)e^{-ikS_{t,\ell}(r)}\right),
\]
with
\[
p_{\ell,t}(r;z)=\left(1+\frac{q_t(r)}{z}-\frac{\ell(\ell+1)}{k^2r^2 }\right)^{1/2},
\]
and
\[
S_{t,\ell}(r)=\int_{r_{\ast}}^rp_{\ell,t}(\varrho;\lambda(t))d\varrho,
\]
and this leads to
\begin{equation}\label{eqn:law_23d}
\frac{\lambda'(t)}{\lambda(t)}=\frac{\displaystyle\int_0^R\frac{q_{\rm rad}(r)-q_0}{p_{\ell,t}(r;\lambda(t))}dr}{\displaystyle\int_0^R\frac{q_t(r)}{p_{\ell,t}(r;\lambda(t))}dr}.
\end{equation}
Similar to Section~\ref{subsec:d=2}, define
\[
H_{\ell}(t)=\int_0^R\left(p_{\ell,t}(r;\lambda(t))-p_{\ell,0}(r;\lambda(0))\right)dr,
\]
and note that, to leading order as $k\rightarrow\infty$,
\begin{align*}
\frac{dH_{\ell}}{dt}(t)&=\frac{1}{2\lambda(t)}\int_0^R\frac{q_{\rm rad}(r)-q_0}{p_{\ell,t}(r;\lambda(t))}dr-\frac{\lambda'(t)}{2\lambda(t)^2}\int_0^R\frac{q_t(r)}{p_{\ell,t}(r;\lambda(t))}dr=0,
\end{align*}
so specifically
\begin{equation}\label{eqn:use2}
\int_0^R\left(p_{\ell,1}(r;\lambda(1))-p_{\ell,0}(r;\lambda(0))\right)dr=H_{\ell}(1)=o(1),\quad k\rightarrow\infty,
\end{equation}
which can be solved for $\lambda(1)$ numerically using a zero-finding algorithm and the precomputed eigenvalue \(\lambda(0)\) of \(k^2V_{q_0}(k)\).

\subsection{Spectral correction} It remains to describe the first-order correction of the computed spectrum for $k^2V_{q_{\rm rad}}(k)$ so that it estimates the spectrum of $k^2V_{q}(k)$. Recall that $q(x)=q_{\rm rad}(|x|)+\eta Q(x)$ with $\|Q\|_{\infty}=1$ and $|\eta|\ll1$, and now let $k^2V_{q_{\rm rad}}(k)\varphi=\lambda\varphi$ and $k^2V_{q}(k)\varphi_{\eta}=\lambda_{\eta}\varphi_{\eta}$, with
\[
\lambda_{\eta}=\lambda+\eta\mu+O(\eta^2),
\]
\[
\varphi_{\eta}=\varphi+\eta\psi+O(\eta^2).
\]
Then $V_{q}(k)=V_{q_{\rm rad}}(k)+\eta V_Q(k)$, and, to first order in $\eta$,
\begin{equation}\label{eqn:abc}
(k^2V_{q_{\rm rad}}(k)-\lambda I)\psi=(\mu I-k^2V_Q(k))\varphi.
\end{equation}
Using again the pairing $\langle u,v\rangle=\int_Du(x)v(x)dx$ and the fact that the eigenpairs of $(k^2V_{q_{\rm rad}}(k))^T$ are $(\lambda,q_{\rm rad}\varphi)$, we pair both sides of~\eqref{eqn:abc} with $q_{\rm rad}\varphi$ and get, to first order in $\eta$, that
\begin{align*}
\mu\int_D q_{\rm rad}(|x|)\varphi(x)^2dx&=k^2\int_Dq_{\rm rad}(|x|)\varphi(x)\int_D\Phi_k^{(d)}(x-y)Q(y)\varphi(y)dydx\\&=k^2\int_DQ(y)\varphi(y)\int_D\Phi_k^{(d)}(x-y)q_{\rm rad}(|x|)\varphi(x)dxdy\\&=\int_DQ(y)\varphi(y)k^2V_{q_{\rm rad}}(k)\varphi(y)dy=\lambda\int_DQ(y)\varphi(y)^2dy.
\end{align*}
Thus, to first order in $\eta$, we have the correction to $\lambda$ given by
\begin{equation}\label{eqn:mu}
\mu=\lambda\frac{\int_DQ(x)\varphi(x)^2dx}{\int_Dq_{\rm rad}(|x|)\varphi(x)^2dx}.
\end{equation}

Recall that our proposed approach is to first solve~\eqref{eqn:use1} (when $d=2$) or~\eqref{eqn:use2} (when $d=3$) to propagate $\lambda(t)$ along the homotopy $q_t=(1-t)q_0+tq_{\rm rad}$, $t\in[0,1]$, and then to correct the found $\lambda(1)$ using $\mu$ from~\eqref{eqn:mu}.

\section{Estimating the optimal preconditioner parameter}\label{sec:estimating}

Let $\widehat\Lambda=\{\widehat\lambda_j(1)\}_{j=1}^N$ be the estimated truncated discrete spectrum of \(k^2V_{q}(k)\). In applications, we will be dealing with finite matrix representations of $k^2V_q(k)$, so it makes sense to assume that $|\lambda|\le\delta$, with known $\delta$, for all elements $\lambda$ of $\sigma(k^2V_q(k))$ that are omitted due to the truncation. For such elements $\lambda$, we have
\[
|1-\gamma+\gamma\lambda|\le|1-\gamma|+\delta|\gamma|,
\]
so, given an admissible set \(\Gamma\subset\C\), we define the optimal preconditioner parameter $\gamma_{\mathrm{opt}}$ by the minimax problem
\begin{equation}\label{eq:gamma_opt_definition}
\gamma_{\mathrm{opt}}\in\arg\min_{\gamma\in\Gamma}
\max_{\widehat\lambda_j\in\widehat\Lambda}
\{|1+\gamma\widehat\lambda_j-\gamma|,|1-\gamma|+\delta|\gamma|\}.
\end{equation}
Note that finding $\gamma_{\rm opt}$ that gives $\min_{\Gamma}\max_{\widehat\Lambda}\{|1+\gamma\widehat\lambda_j-\gamma|,|1-\gamma|+\delta|\gamma|\}<1$ does not necessarily guarantee that the spectral radius of $M_{\gamma_{\rm opt}}$ is less than one, since for any true eigenvalue $\lambda_j$ of $k^2V_q(k)$, approximated by $\widehat\lambda_j\in\widehat\Lambda$, we have
\[
|1-\gamma+\gamma\lambda_j|\le|1-\gamma+\gamma\widehat\lambda_j|+|\gamma||\lambda_j-\widehat\lambda_j|.
\]
However, in our numerical tests the solution of the minimax problem sufficed to make the modified Born series convergent. In the computations in Section~\ref{sec:numerical}, we take \(\Gamma\) to be a rectangular box in the first quadrant,
\[
\Gamma=\bigl\{\gamma=u+iv:\ 0\le u\le \gamma_{\Re}^{\max},\ 0\le v\le \gamma_{\Im}^{\max}\bigr\}.
\]
For the 'exact' preconditioner parameter used in Section~\ref{sec:numerical}, the optimization problem \eqref{eq:gamma_opt_definition} is solved by Sequential Least Squares Programming (SLSQP) optimization with the estimated eigenvalues replaced by the numerically computed reference eigenvalues. Next, introducing the auxiliary variable \(r\), one may rewrite \eqref{eq:gamma_opt_definition} as $\min_{\gamma\in\Gamma, r\ge 0} r$
subject to
\[
\max\{|1+\gamma(\widehat\lambda_j(1)-1)|,|1-\gamma|+\delta|\gamma|\}\le r,\qquad j=1,\dots,N.
\]
Writing \(\gamma=u+iv\) and \(\mu_j=\widehat\lambda_j(1)-1=\alpha_j+i\beta_j\), each constraint becomes
\[
\left\|\begin{pmatrix}
    \Re(1+\gamma\mu_j)\\\Im(1+\gamma\mu_j)
\end{pmatrix}\right\|_2=\left\|\begin{pmatrix}1+\alpha_j u-\beta_j v\\\beta_j u+\alpha_j v\end{pmatrix}\right\|_2\le r,
\]
\[
\left\|\begin{pmatrix}
    1-u\\v
\end{pmatrix}\right\|_2\le s_0,\quad \left\|\begin{pmatrix}
    u\\v
\end{pmatrix}\right\|_2\le s_1,\quad s_0+\delta s_1\le r.
\]
This is a convex second-order cone program in the real variables \((u,v,r,s_0,s_1)\). Given any candidate parameter \(\gamma\), we evaluate its performance on the estimated spectrum through the contraction metric $\varrho(\gamma)
:=\max_{\widehat\lambda_j\in\widehat\Lambda}\{|1+\gamma\widehat\lambda_j-\gamma|,|1-\gamma|+\delta|\gamma|\}$. The optimal value with the selected parameter is $\varrho_{\mathrm{opt}}=\varrho\bigl(\gamma_{\mathrm{opt}}\bigr)$. This metric estimates the quality of the preconditioner, in that a smaller value of this objective predicts a better worst-case asymptotic convergence factor.

A second, more interpretable metric is the relative residual
\begin{equation}\label{eqn:rr}
r(\gamma,N):=\frac{\|M_{\gamma}(k)^{N+1}V_1f\|_2}{\|V_1f\|_2}
\end{equation}
associated with the truncated solution
\[
u_N:=\sum_{n=0}^NM_{\gamma}(k)^n\gamma V_1f
\]
of the equation $(I-k^2V_q(k))u=V_1f$. Indeed, since $I-M_{\gamma}(k)=\gamma(I-k^2V_q(k))$, we have
\[
\gamma(I-k^2V_q(k))u_N=\gamma(I-M_{\gamma}(k)^{N+1})V_1f,
\]
so $M_{\gamma}(k)^{N+1}V_1f=-((I-k^2V_q(k))u_N-V_1f)$ is the residual associated with $u_N$. Note that $r(\gamma,N)$ is pre-computable for any $\gamma$ and $N$, and we can use it for a sum truncation criterion: given $\gamma$ and a tolerance $\varepsilon>0$, the stopping index is the first $N$ satisfying $r(\gamma,N)<\varepsilon$. Numerical experiments show that, for fixed $\varepsilon$, the choice of $\gamma$ greatly affects the value of the stopping index. Moreover, we note that the spectral radius of the operator $M_{\gamma}(k)=(1-\gamma)I+\gamma k^2V_q(k)$ approaches 1 as $\gamma\rightarrow0$, so we expect that the stopping index approaches infinity as $\gamma\rightarrow0$ and $\varepsilon$ is kept fixed. On the other hand, the spectral radius of $M_{\gamma}(k)$ being less than one requires a sufficiently small $|\gamma|$. Thus we expect that, given a tolerance $\varepsilon$, there will in general be an optimum overrelaxation $\gamma$ resulting in a smallest stopping index and hence the fastest MBS solution of the equation $(I-k^2V_q(k))u=V_1(k)f$. Having fixed $\gamma$, we can now either calculate an estimate $N_{\rm estimate}$ of the stopping index and then compute the MBS in~\eqref{eqn:MBS} with $N_{\rm estimate}+1$ summands, or we can test the stopping criterion $r(\gamma,N)<\varepsilon$ more or less frequently as the MBS is being computed.

\section{Numerical examples}\label{sec:numerical}
\begin{figure}
\centering
\includegraphics[width=0.49\textwidth]{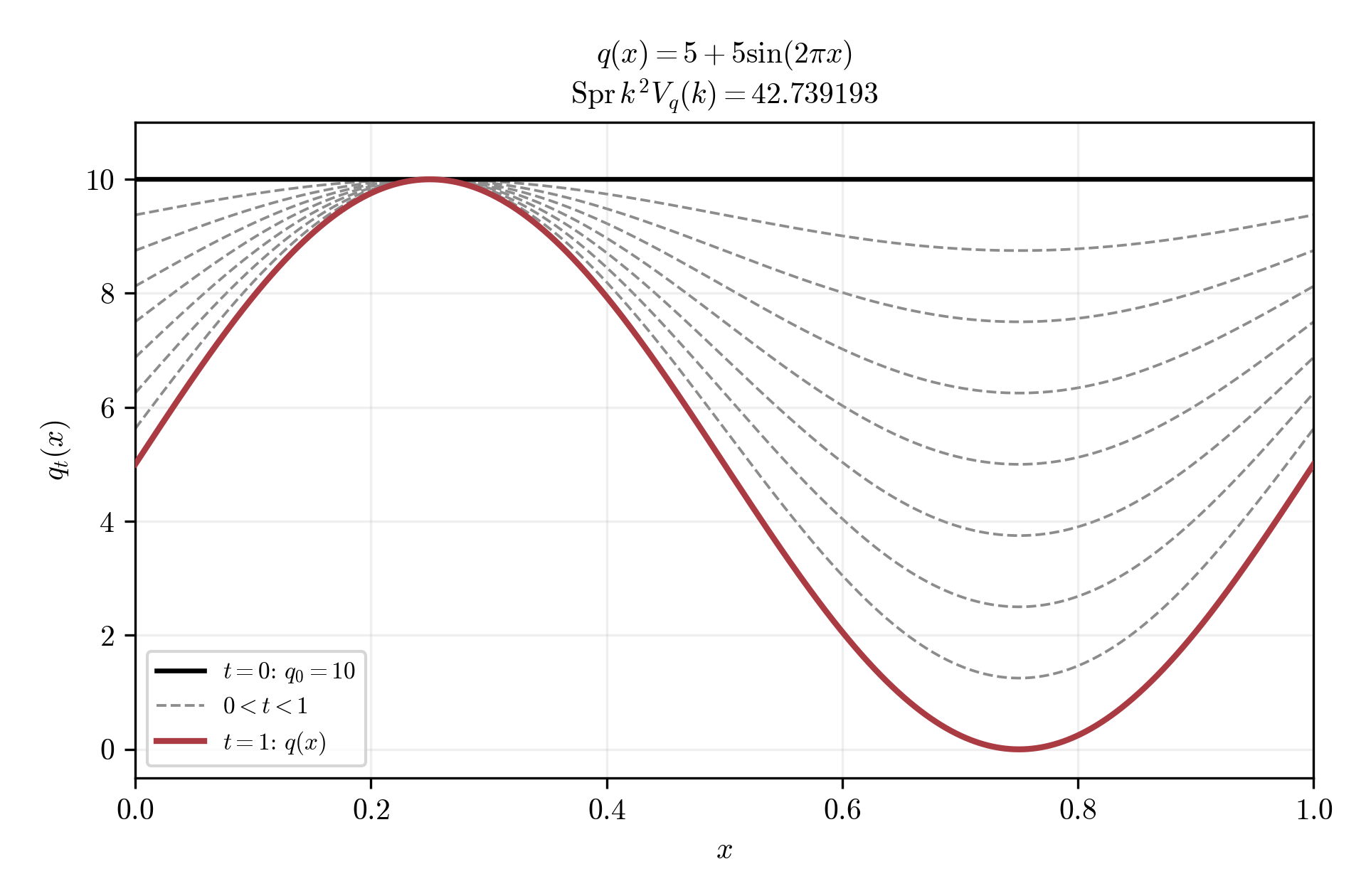}\\
\includegraphics[width=0.99\textwidth]{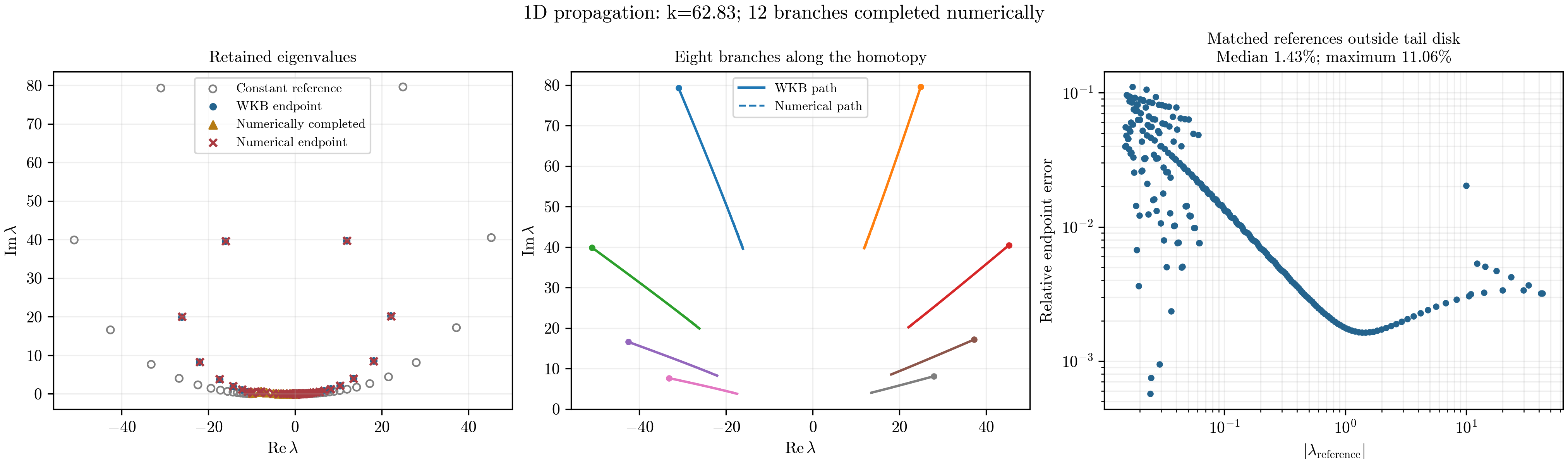}\\
\includegraphics[width=0.99\textwidth]{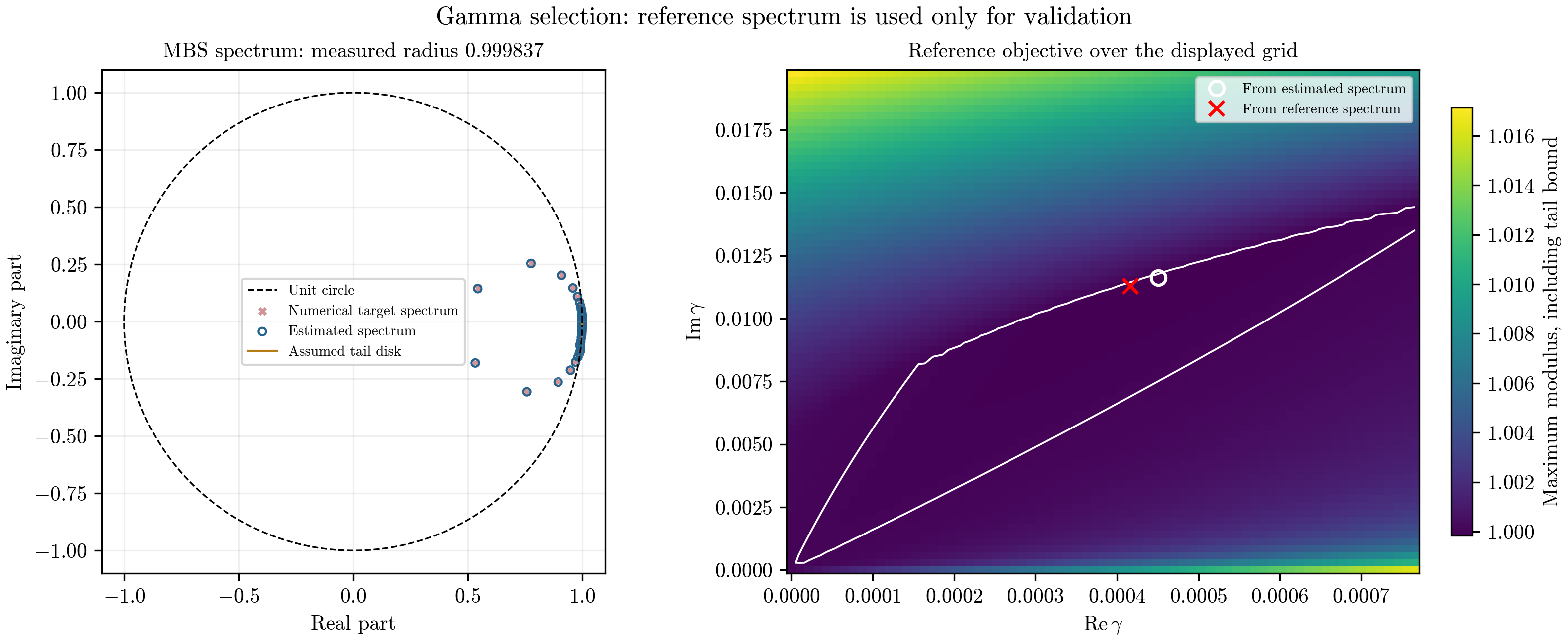}
\caption{The homotopy $q_t$, spectral propagation, and selection of preconditioner $\gamma$ for the case $d=1$, $q(x)=5+5\sin(2\pi x)$. The spectral radius of $k^2V_q(k)$ is approximately 42.7. The white contour is
$\max\{\max_j|1-\gamma+\gamma\lambda_j^{\rm ref}|,
|1-\gamma|+\delta|\gamma|\}=1$,
where $\lambda_j^{\rm ref}$ are the numerically computed eigenvalues
of $k^2V_q(k)$. The expression is below one inside the contour.}\label{fig:1D_sine_1}
\end{figure}

\begin{figure}
\centering
\includegraphics[width=0.99\textwidth]{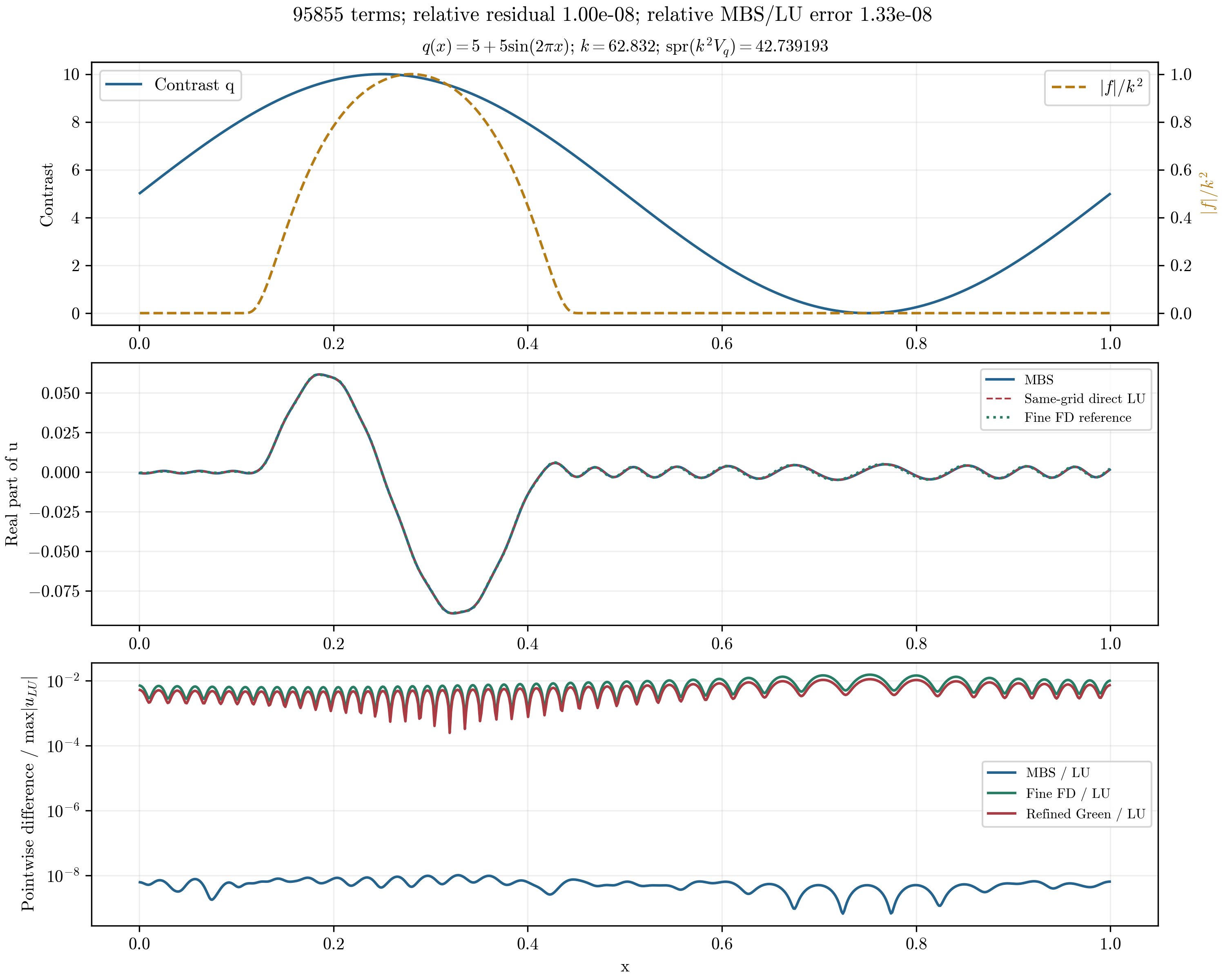}\\
\includegraphics[width=0.99\textwidth]{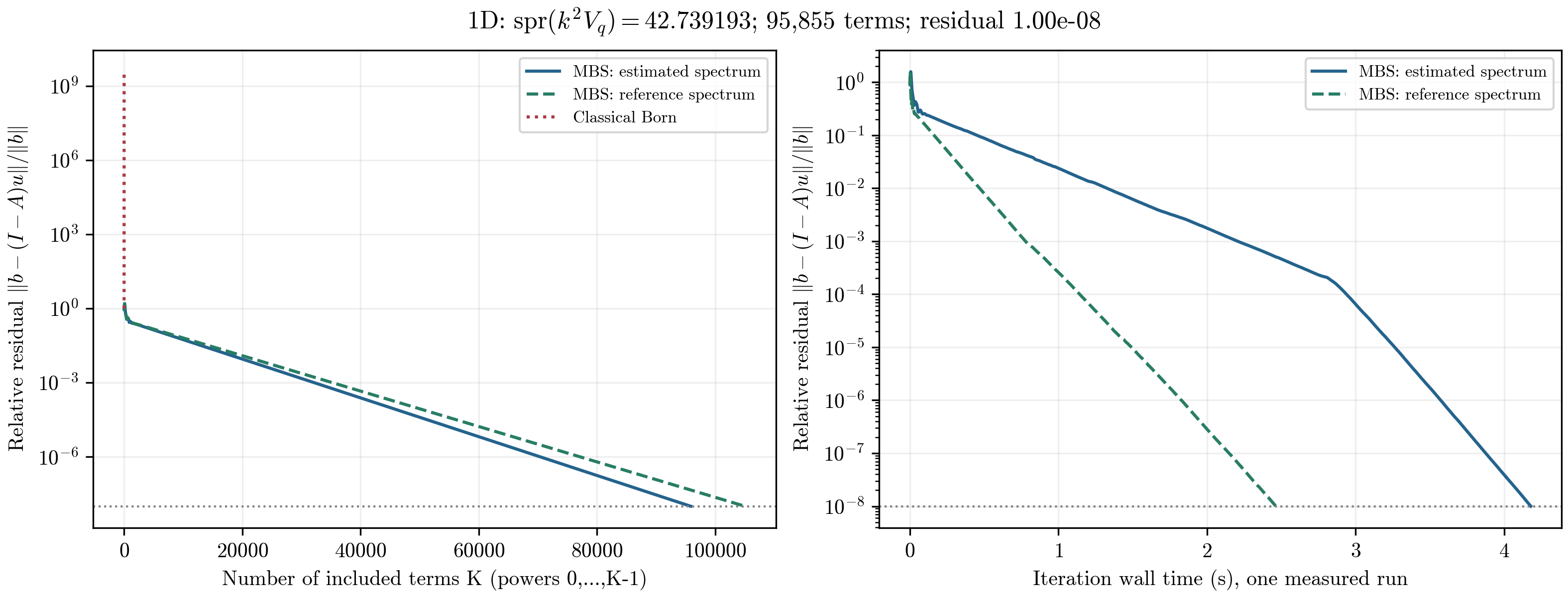}
\caption{The considered medium function $q$ and source $f$; comparison of MBS, direct LU, and FD solutions of~\eqref{eqn:PDE_d1}; relative errors in the numerical solutions; relative residual in MBS solution.}\label{fig:1D_sine_2}

\end{figure}

\begin{figure}
\centering
\includegraphics[width=0.49\textwidth]{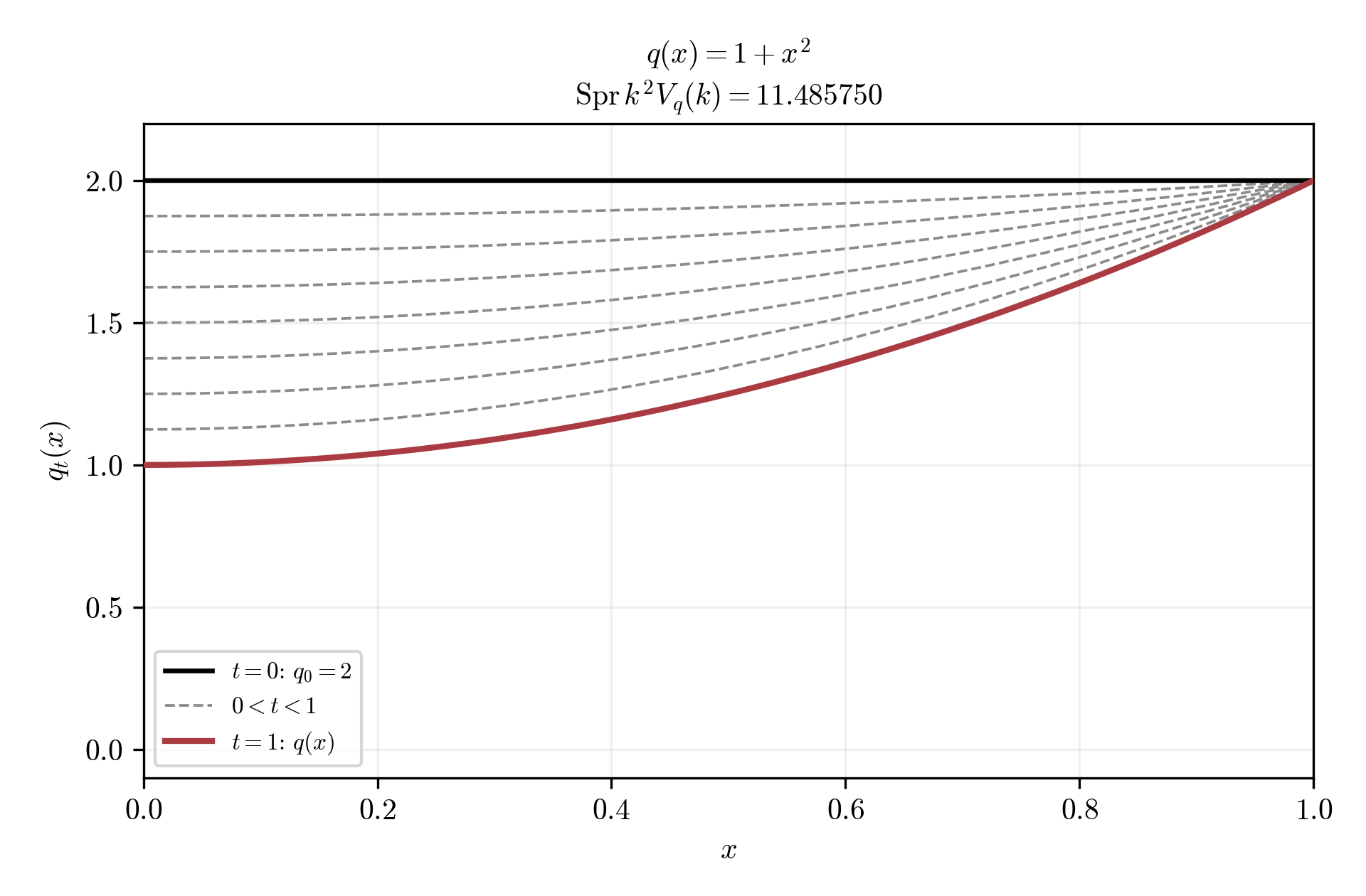}
\includegraphics[width=0.99\textwidth]{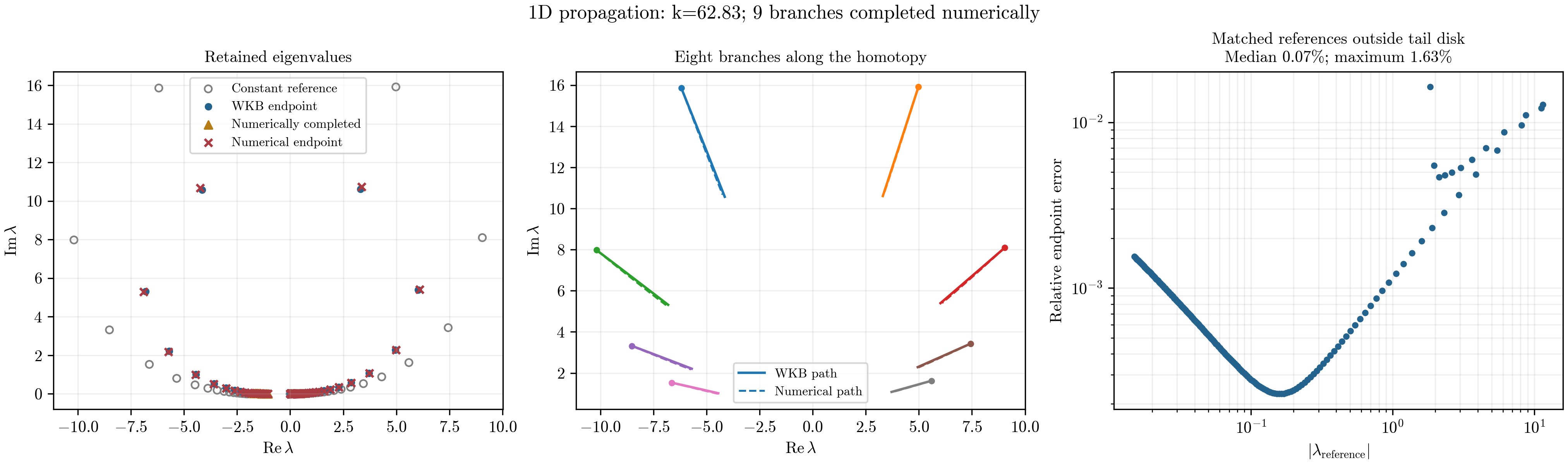}\\
\includegraphics[width=0.99\textwidth]{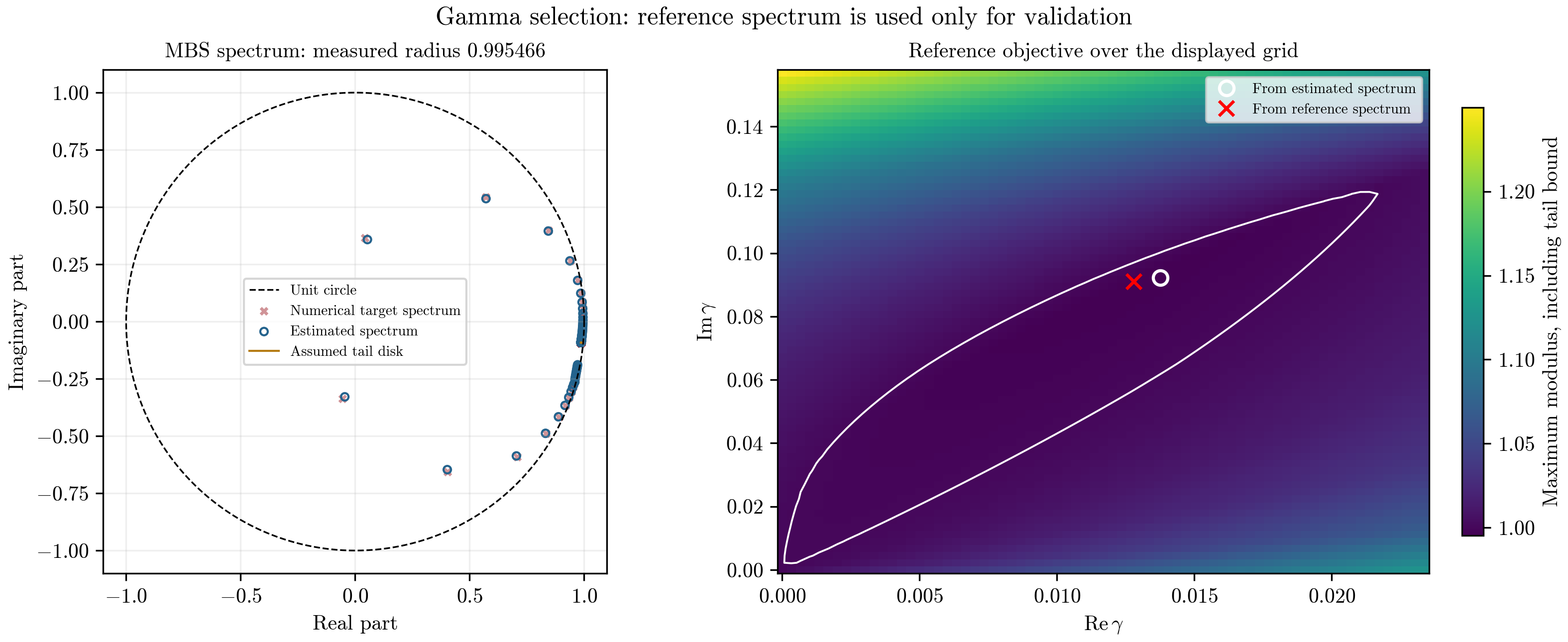}
\caption{The homotopy $q_t$, spectral propagation, and selection of preconditioner $\gamma$ for the case $d=1$, $q(x)=1+x^2$. The spectral radius of $k^2V_q(k)$ is approximately 11.5. The white contour is
$\max\{\max_j|1-\gamma+\gamma\lambda_j^{\rm ref}|,
|1-\gamma|+\delta|\gamma|\}=1$,
where $\lambda_j^{\rm ref}$ are the numerically computed eigenvalues
of $k^2V_q(k)$. The expression is below one inside the contour.}\label{fig:1D_square_1}
\end{figure}

\begin{figure}
\centering
\includegraphics[width=0.99\textwidth]{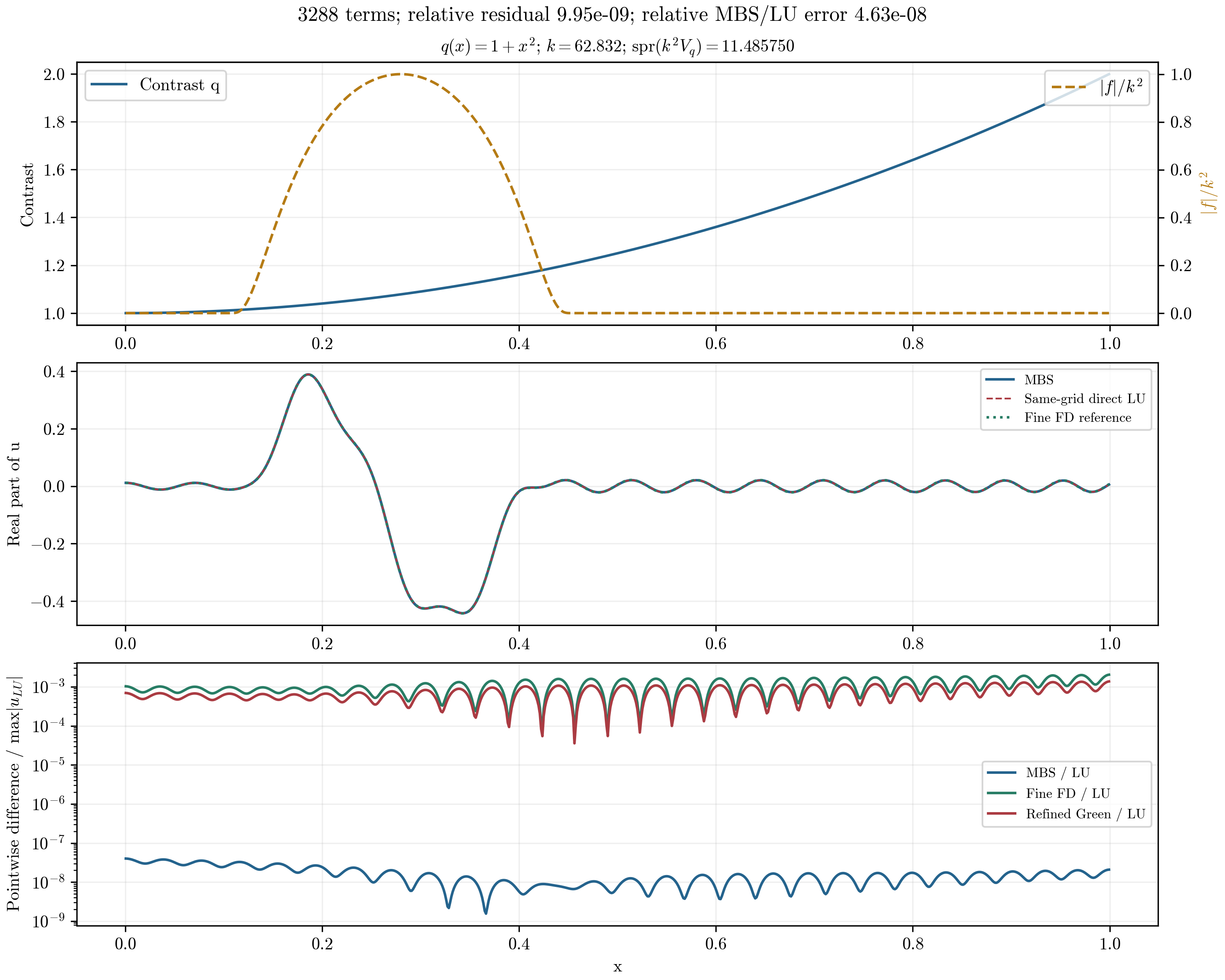}\\
\includegraphics[width=0.99\textwidth]{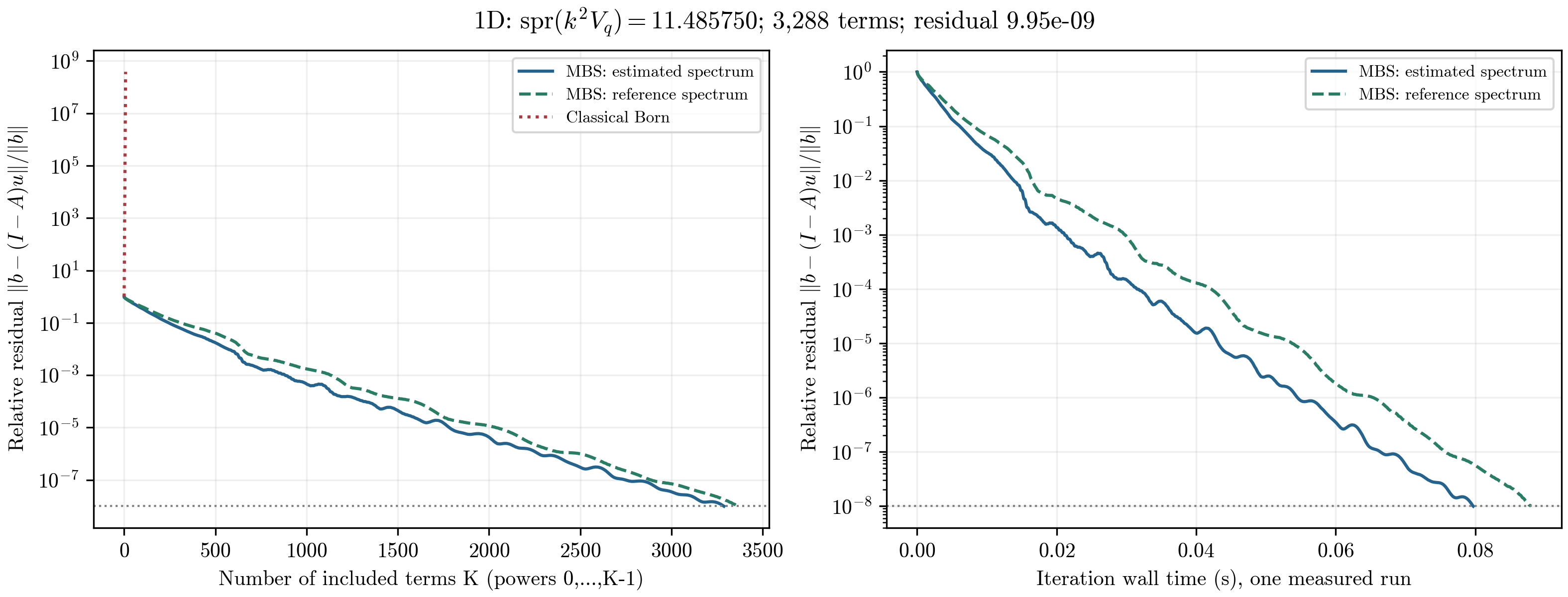}
\caption{The considered medium function $q$ and source $f$; comparison of MBS, direct LU, and FD solutions of~\eqref{eqn:PDE_d1}; relative errors in the numerical solutions; relative residual in MBS solution.}\label{fig:1D_square_2}
\end{figure}

\begin{figure}
\centering
\includegraphics[width=0.49\textwidth]{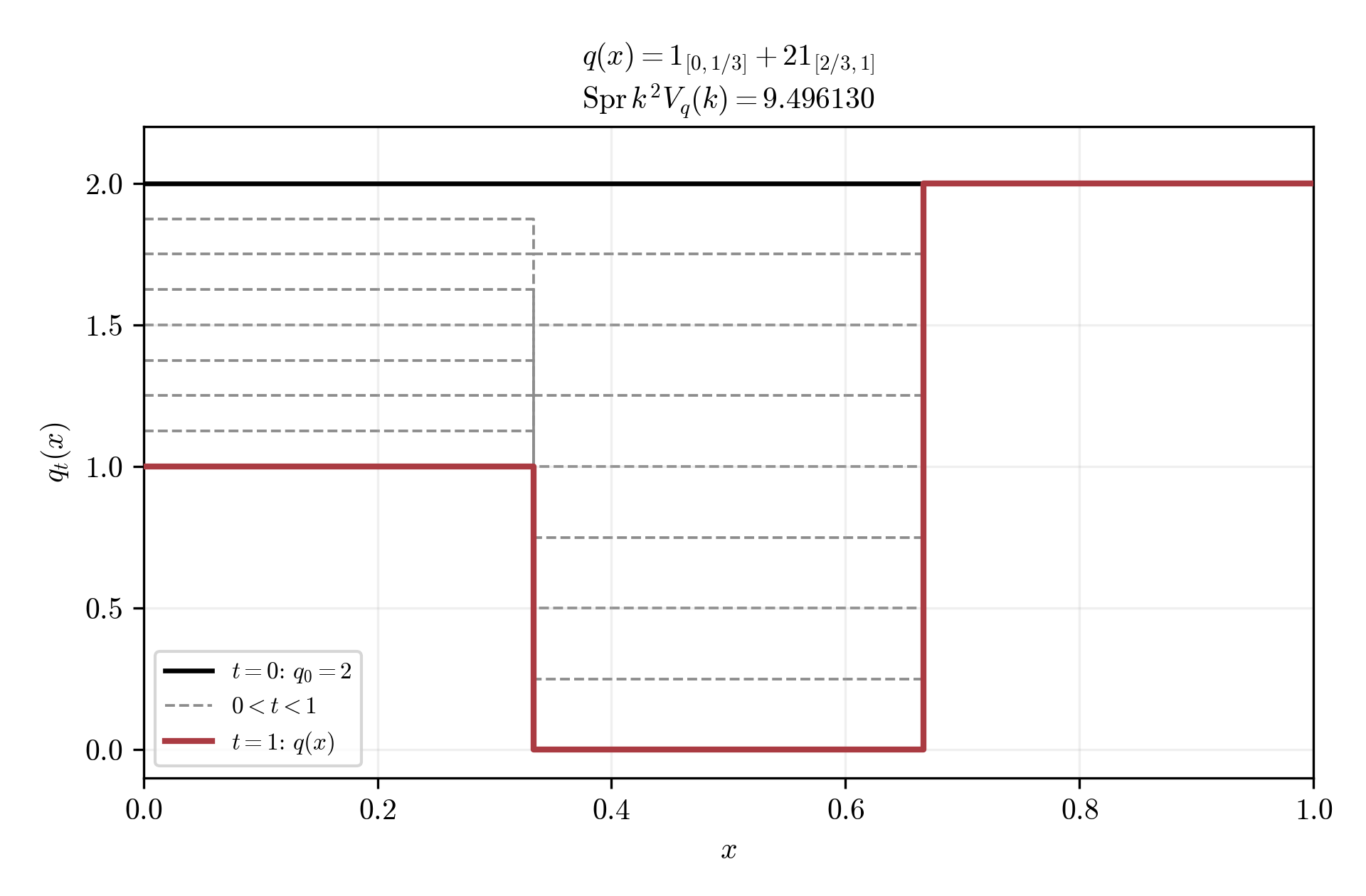}\\
\includegraphics[width=0.99\textwidth]{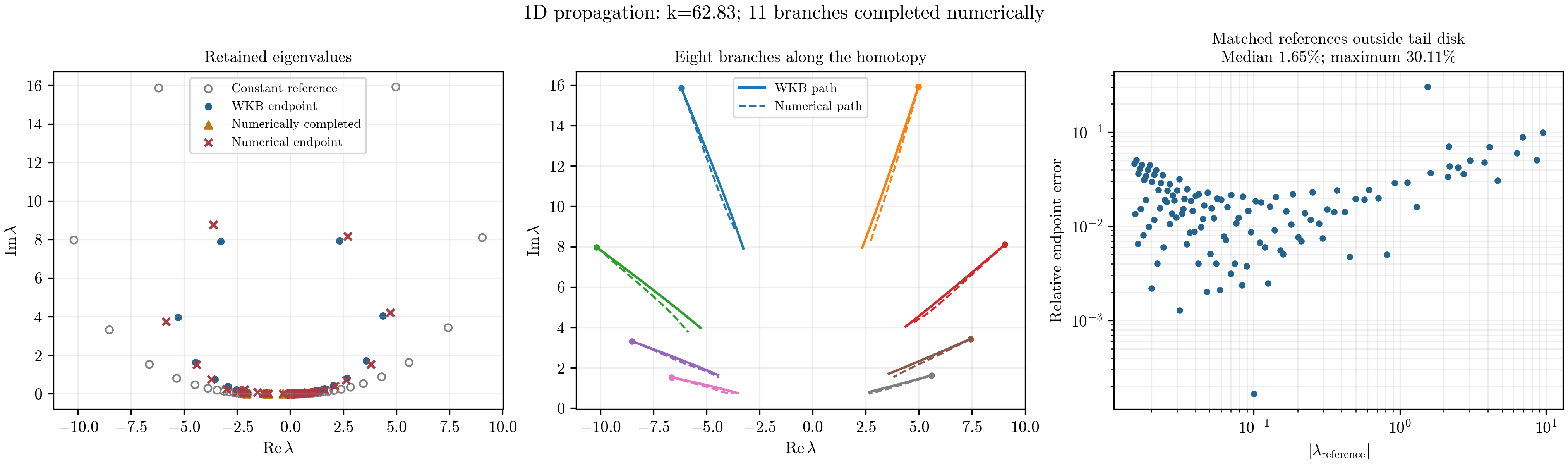}\\
\includegraphics[width=0.99\textwidth]{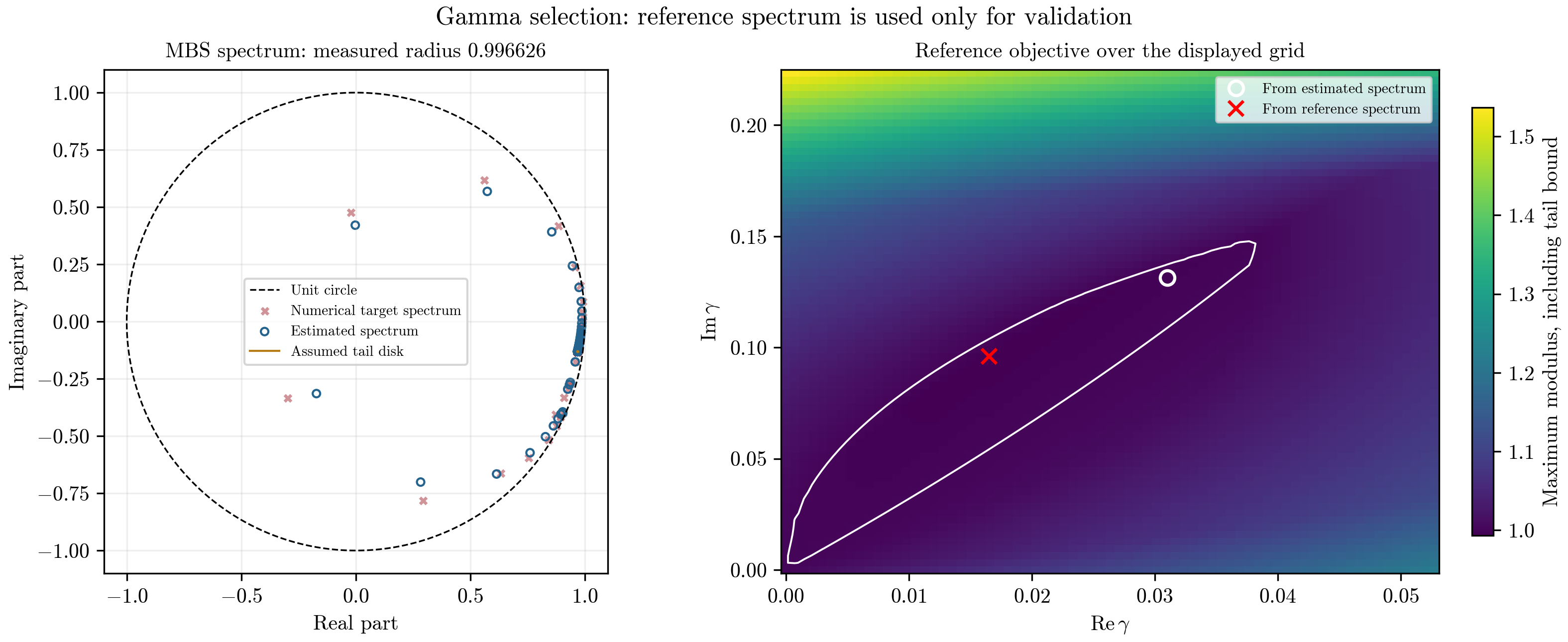}
\caption{The homotopy $q_t$, spectral propagation, and selection of preconditioner $\gamma$ for the case $d=1$, $q(x)=\bm1_{[0,1/3]}+2\cdot\bm1_{[2/3,1]}$. The spectral radius of $k^2V_q(k)$ is approximately 9.5. The white contour is
$\max\{\max_j|1-\gamma+\gamma\lambda_j^{\rm ref}|,
|1-\gamma|+\delta|\gamma|\}=1$,
where $\lambda_j^{\rm ref}$ are the numerically computed eigenvalues
of $k^2V_q(k)$. The expression is below one inside the contour.}\label{fig:1D_steps_1}
\end{figure}

\begin{figure}
\centering
\includegraphics[width=0.99\textwidth]{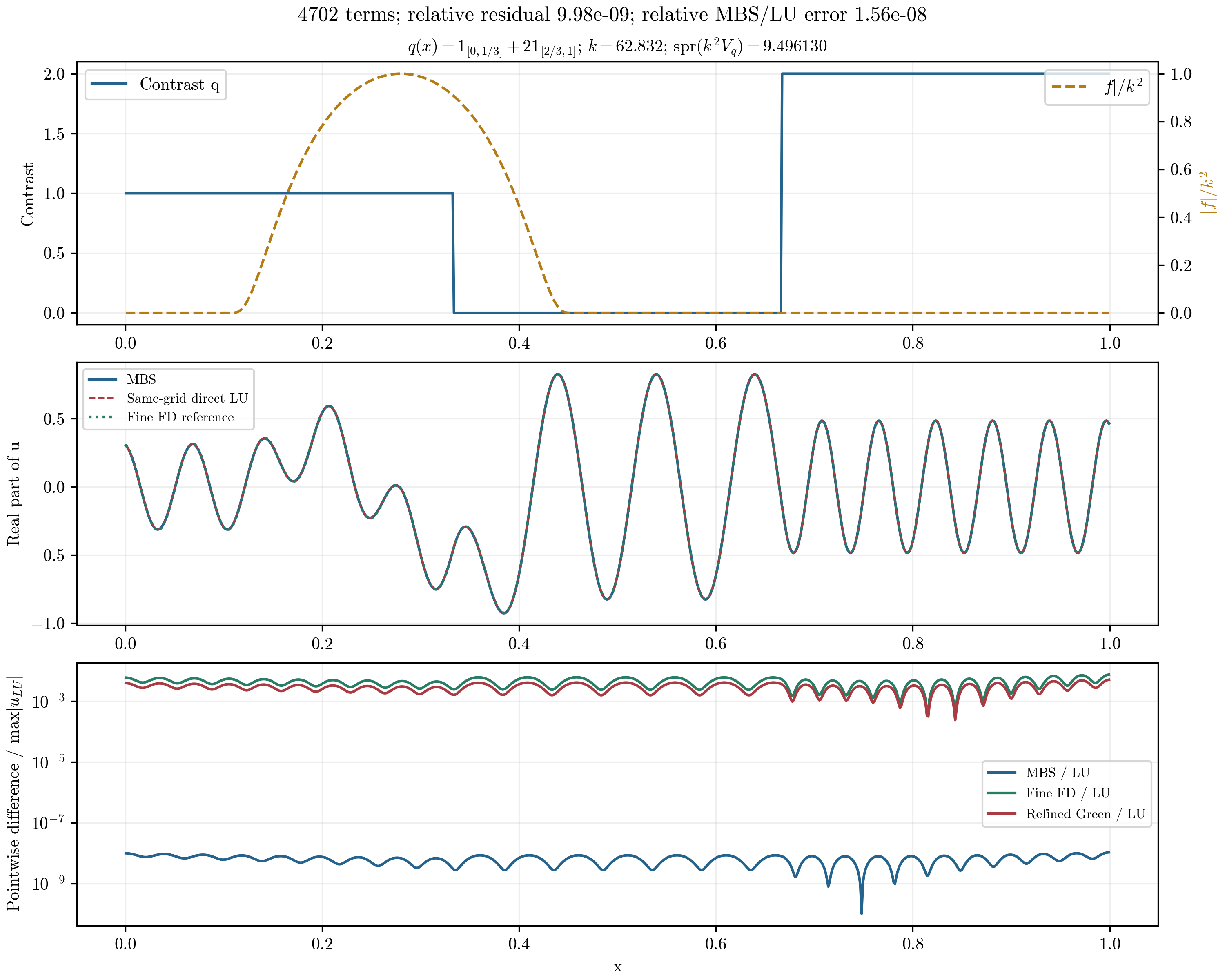}\\
\includegraphics[width=0.99\textwidth]{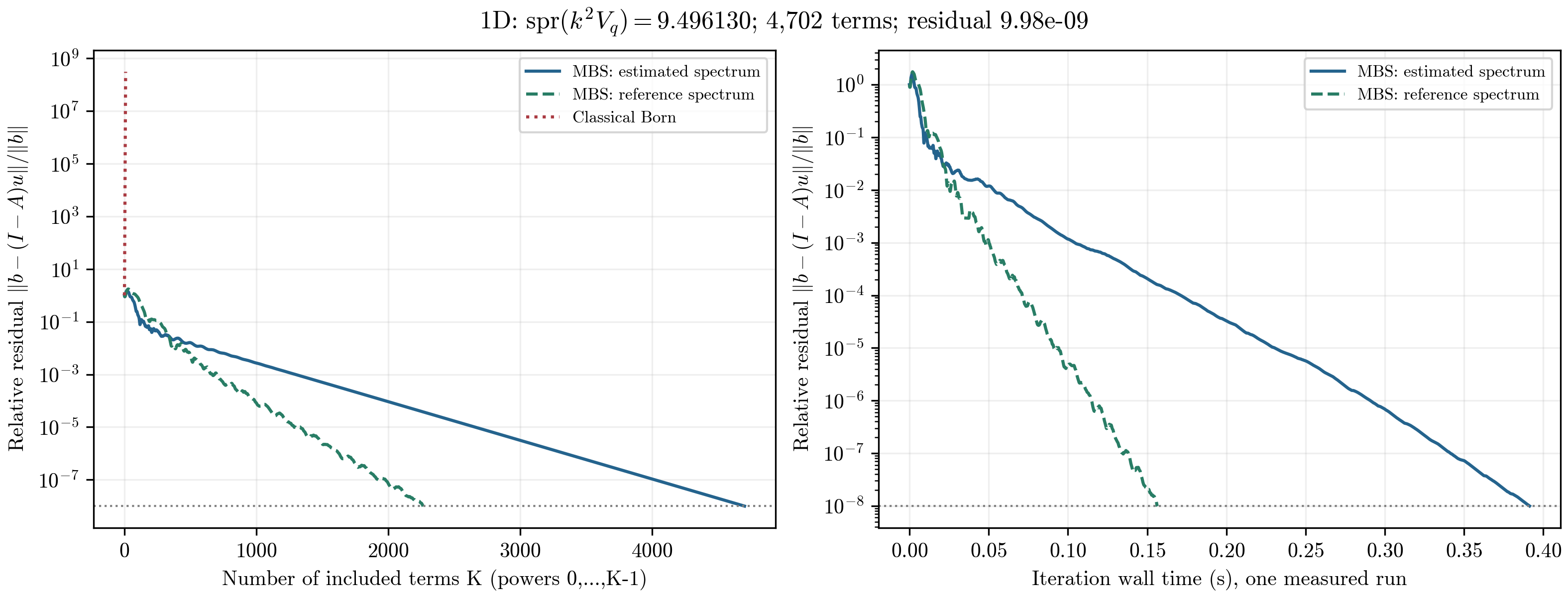}
\caption{The considered medium function $q$ and source $f$; comparison of MBS, direct LU, and FD solutions of~\eqref{eqn:PDE_d1}; relative errors in the numerical solutions; relative residual in MBS solution.}\label{fig:1D_steps_2}
\end{figure}

We next illustrate the proposed accelerated modified Born series approach using numerical examples in one and two dimensions. In one dimension, we consider the media $q(x)=5+5\sin(2\pi x)$, $q(x)=1+x^2$, and $q(x)=1_{[0,1/3]}+2\cdot1_{[2/3,1]}$, $x\in[0,1]$, with the corresponding spectral radii of $k^2V_q(k)$ being approximately $42.7$, $11.5$, and $9.5$, respectively, and with the same source used for all three
contrasts,
\[
f(x)=k^2\beta\!\left(\frac{|x-0.28|}{0.18}\right)e^{0.3ikx},
\qquad k=20\pi,
\]
where
\[
\beta(s)=
\begin{cases}
\exp\!\left(1-\dfrac{1}{1-s^2}\right), & 0\le s<1,\\
0, & s\ge1.
\end{cases}
\]
Figures~\ref{fig:1D_sine_1}--\ref{fig:1D_steps_2} demonstrate the applicability of the approach in one dimension for these cases where the spectral radius of $k^2V_q(k)$ is significantly larger than one.

\begin{figure}
    \centering
    \includegraphics[width=0.99\linewidth]{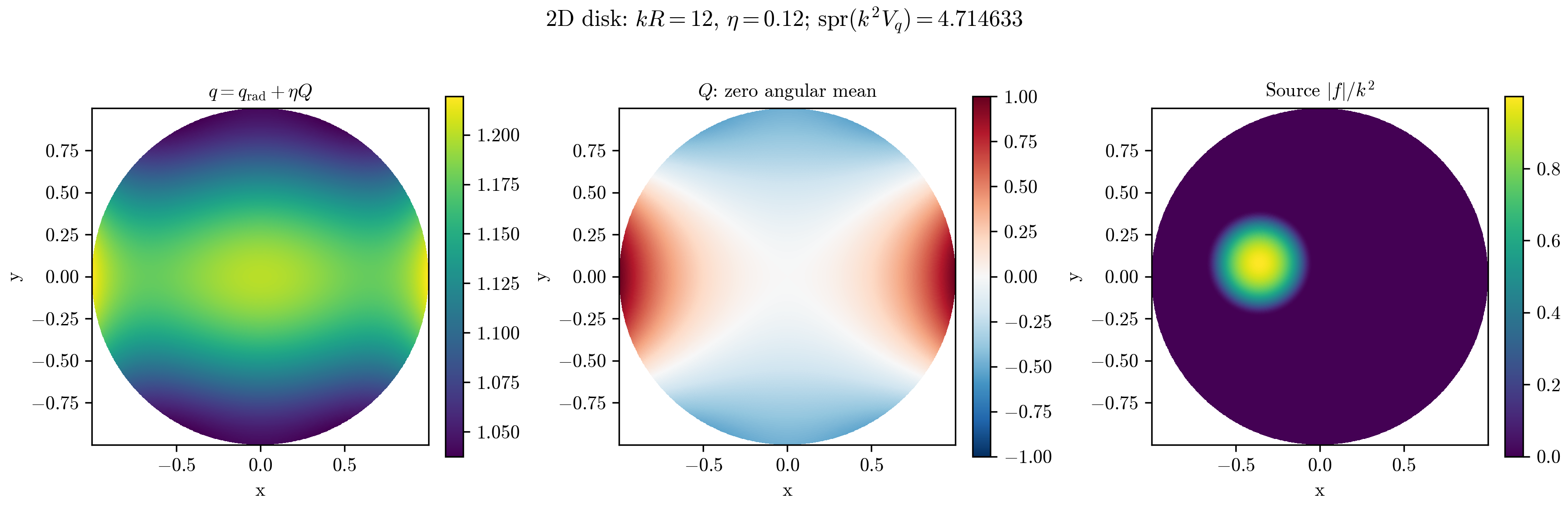}\\
    \includegraphics[width=0.99\linewidth]{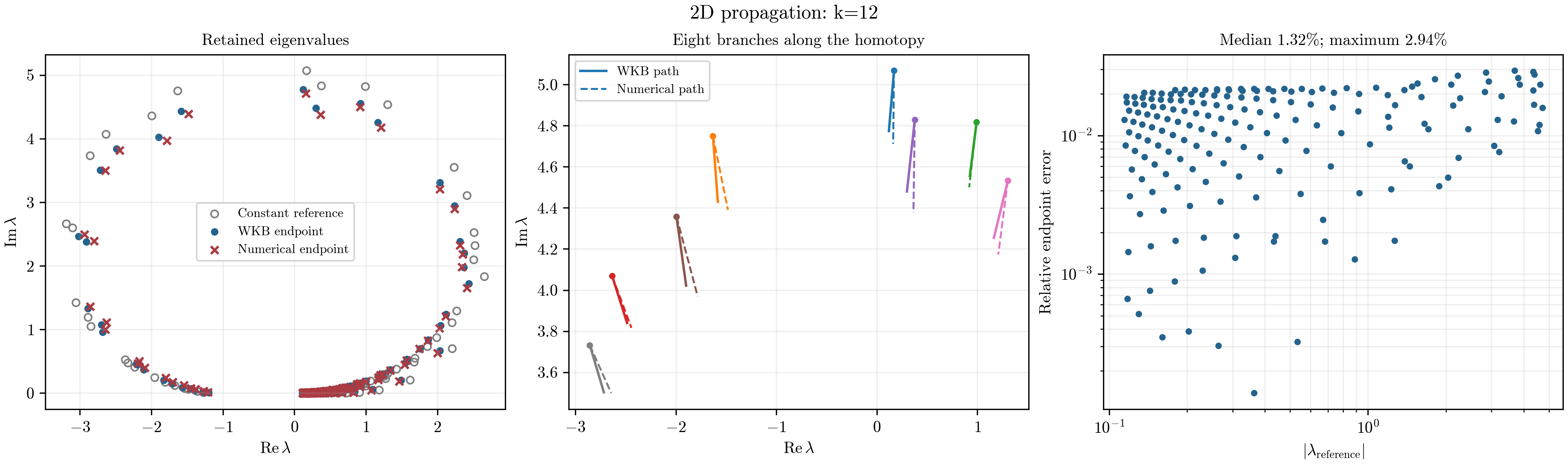}\\
    \includegraphics[width=0.99\linewidth]{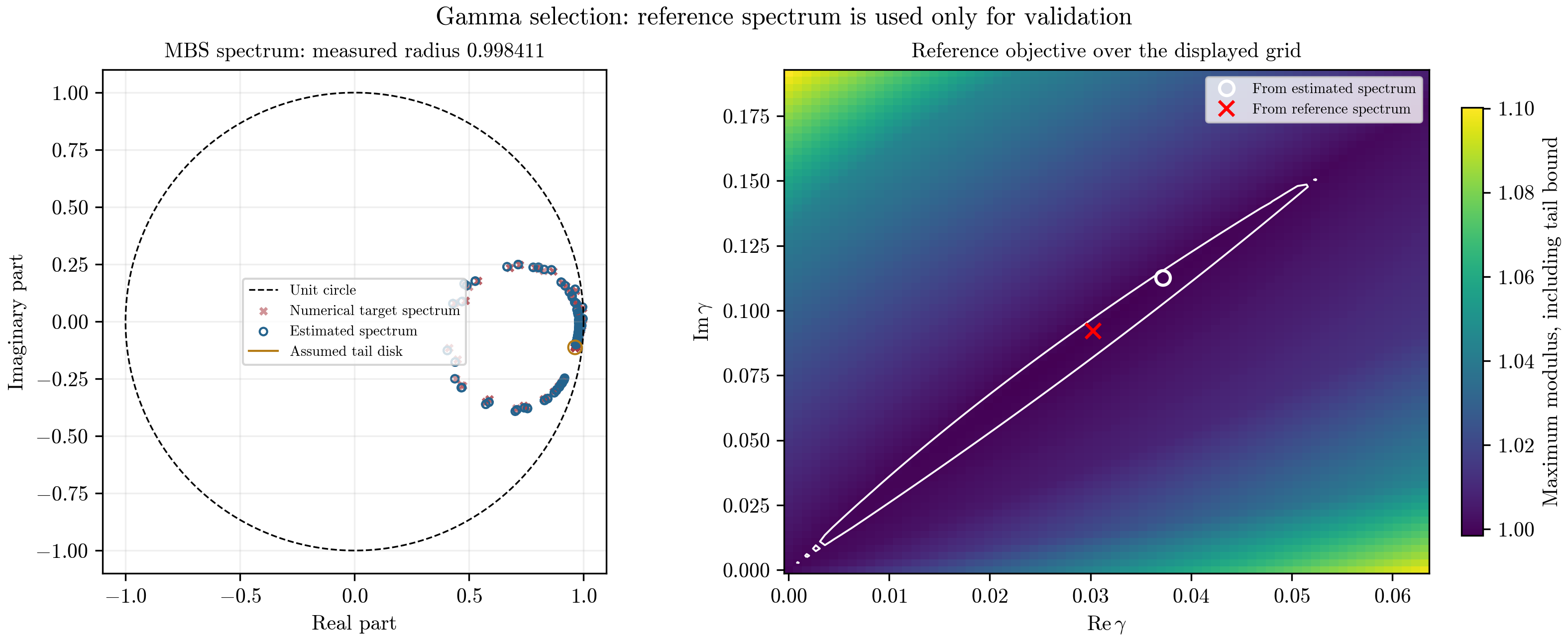}
    \caption{The considered weakly nonradial medium $q$ and source $f$ for the example of case $d=2$; spectral propagation; selection of preconditioner $\gamma$. The spectral radius of $k^2V_q(k)$ is approximately 4.7. The white contour is
$\max\{\max_j|1-\gamma+\gamma\lambda_j^{\rm ref}|,
|1-\gamma|+\delta|\gamma|\}=1$,
where $\lambda_j^{\rm ref}$ are the numerically computed eigenvalues
of $k^2V_q(k)$. The expression is below one inside the contour.}\label{fig:2D_1}
\end{figure}

\begin{figure}
    \centering
    \includegraphics[width=0.99\linewidth]{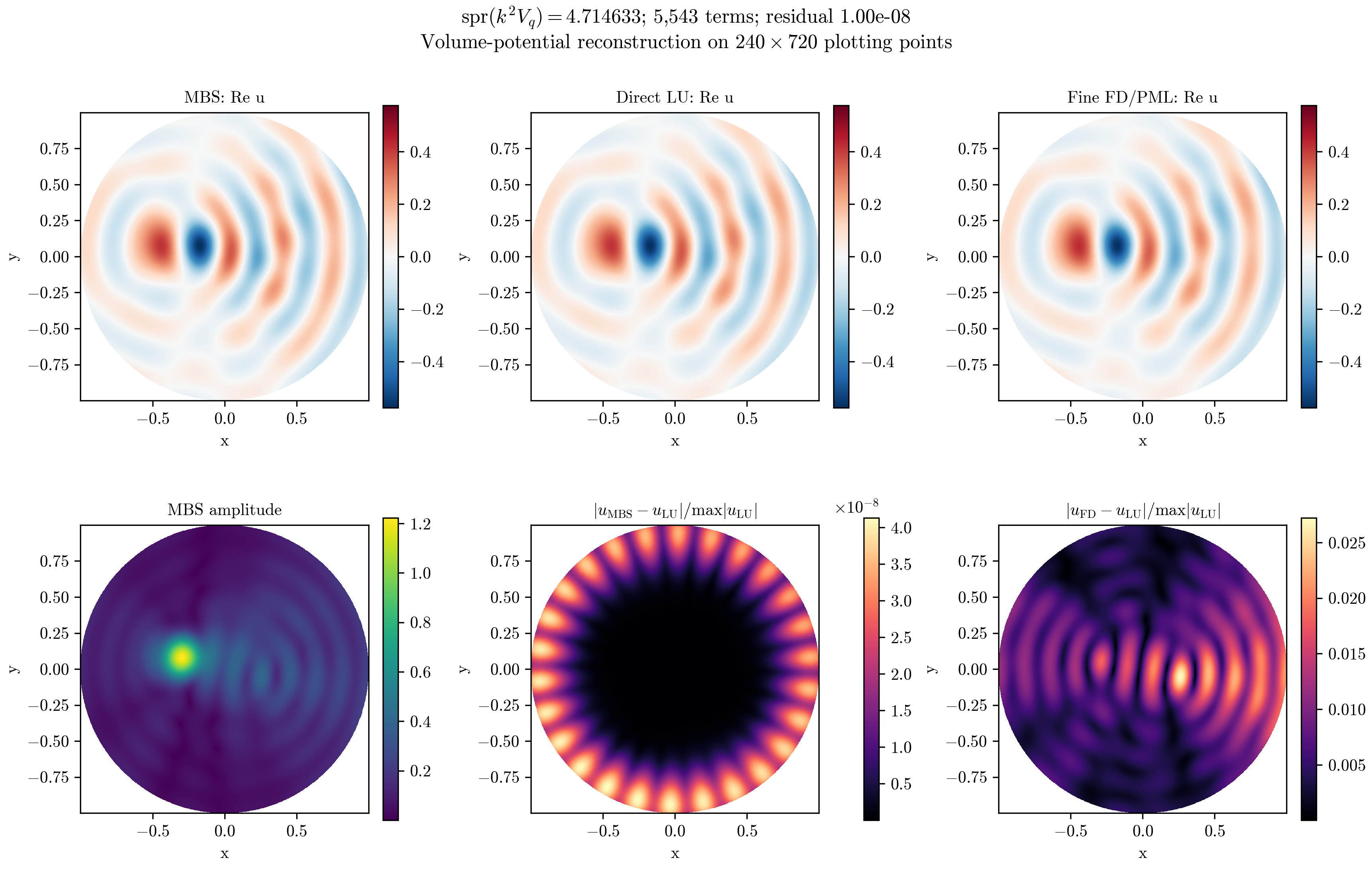}\\
    \includegraphics[width=0.99\linewidth]{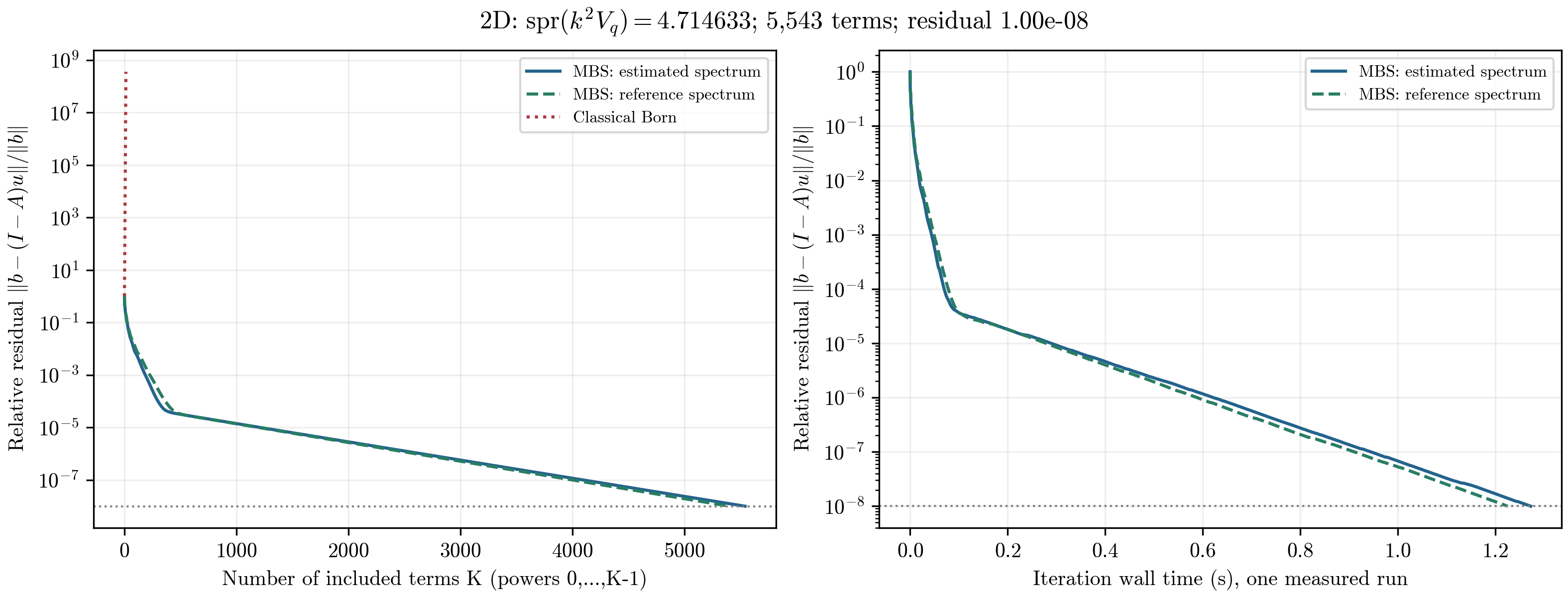}
    \caption{Comparison of MBS, direct LU, and FD solutions of~\eqref{eqn:PDE_d23}; relative errors; relative residual in MBS solution.}\label{fig:2D_2}
\end{figure}

Figures~\ref{fig:2D_1}--\ref{fig:2D_2} support the use of the approach in two dimensions, where we have the medium
\[
q(x)=
\begin{cases}
1.1+0.1(1-|x|^2)^2
+\dfrac{0.12}{1.35}
\left[x_1^2-x_2^2
+0.35\left(x_1^4-6x_1^2x_2^2+x_2^4\right)\right],
& |x|\le1,\\[2mm]
0, & |x|>1,
\end{cases}
\]
and specifically $\eta=0.12$ and $\|Q\|_{\infty}=1$, as well as the source
\[
f(x_1,x_2)=k^2
\beta\!\left(
\frac{\sqrt{(x_1+0.36)^2+(x_2-0.08)^2}}{0.34}
\right)e^{0.4ikx_1},
\qquad k=12,
\]
supported in the disk of radius $0.34$ centred at
$(-0.36,0.08)$.

\begin{figure}
    \centering
    \includegraphics[width=0.99\linewidth]{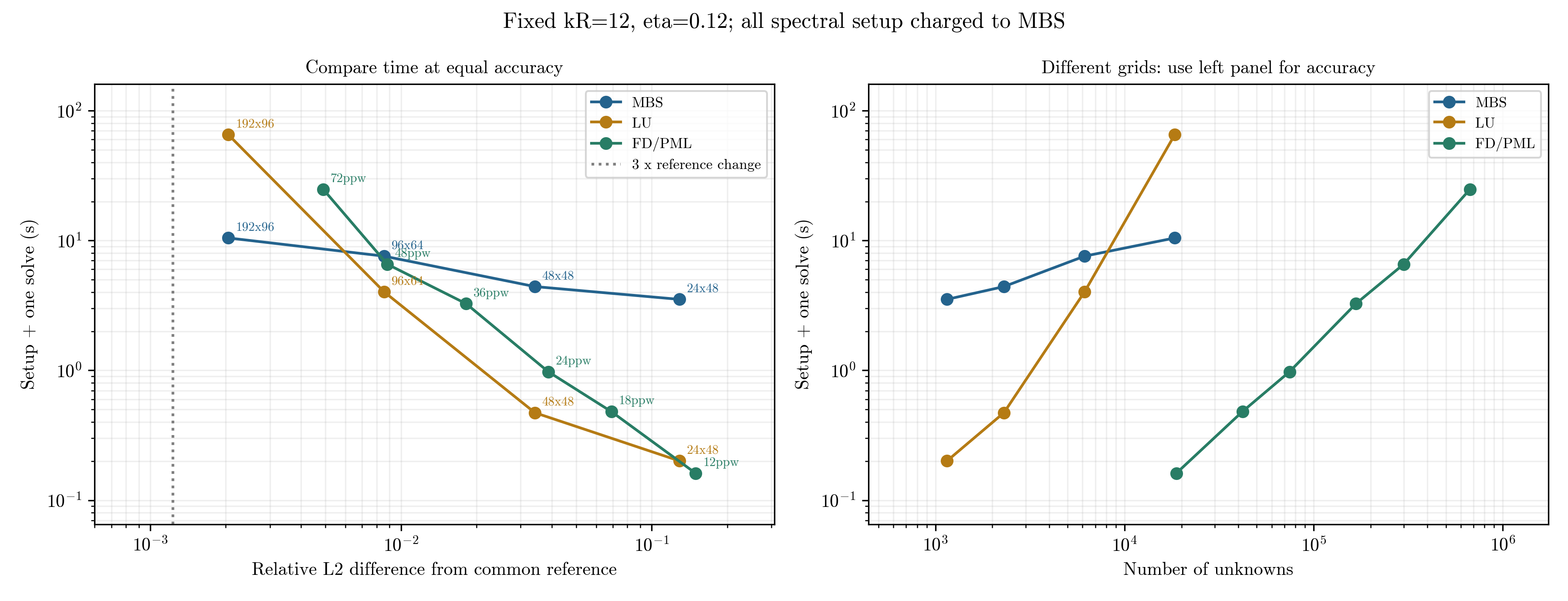}
    \caption{Execution times for the numerical methods applied to the given two-dimensional medium problem. Our accelerated MBS outperforms both direct LU and FD at sufficient accuracy.}\label{fig:2D_3}
\end{figure}

At the more demanding accuracy levels tested in Figure~\ref{fig:2D_3}, MBS achieves a lower setup-plus-solve time than the two implemented competing methods.

The experiments use $k=20\pi$ and $768$ uniform cells in
one dimension, and $k=12$ on the unit disk with a
$48\times48$ radial--angular grid in two dimensions.
MBS and direct LU use the same volume-potential
discretization, with piecewise constant radial elements
and, in two dimensions, angular Fourier collocation.
Fast matrix applications use cumulative sums and FFTs.
The independent competitor uses second-order finite
differences, with outgoing Robin conditions in one
dimension and a perfectly matched layer in two dimensions.

Spectral propagation uses eight homotopy steps. When Newton iteration for a one-dimensional WKB action equation
fails to reach the residual tolerance $2\cdot10^{-8}$, we restart
from the initial eigenpair and follow the discretized eigenvalue
problem using shifted inverse iteration with adaptive homotopy
steps. The resulting endpoint replaces the WKB estimate when
selecting $\gamma$. This additional computation is included
in the MBS setup time. Initial eigenvalues are retained above
modulus thresholds $0.01$ and $0.12$ in one and two
dimensions, respectively; the corresponding assumed
tail radii are $0.015$ and $0.25$. The parameter $\gamma$
is selected by Sequential Least Squares Programming (SLSQP), with an additional spectral-error
margin $0.015|\gamma|$ in two dimensions.
MBS stops at relative residual $10^{-8}$.

The refinement study uses polar grids
$(24,48)$, $(48,48)$, $(96,64)$, and $(192,96)$,
and finite-difference grids with $12$--$72$ points per
shortest estimated wavelength. Errors are relative
$L^2$ differences, including the polar area weight,
from a common MBS reference on a $768\times96$ grid
with residual below $10^{-10}$. Additional
$384\times96$ and $384\times80$ calculations assess
reference sensitivity. Agreement with LU on the same
grid measures iteration error rather than discretization
accuracy.

Timings include setup and one solve, charging the full
spectral setup to every MBS result. Factorization and
solve times are medians of three repetitions with one
BLAS/OpenMP thread.


\section{Conclusion}~\label{sec:conclusion}
We have developed an accelerated Born series method for the numerical solution of scattering problems governed by the variable-coefficient Helmholtz equation with nonnegative-valued, compactly supported, bounded, weakly non-radial media. We also provided a practical numerical procedure for selecting the overrelaxation parameter. The method is effective in the high-frequency regime, for both weak and strong scattering, and has been successfully validated against direct integral-equation and finite-difference computations.

\begin{acknowledgements} This work was funded by Villum Foundation research grant no. 58857.
\end{acknowledgements}

\bibliographystyle{siam}
\bibliography{references}
\end{document}